\documentclass[12pt,reqno]{amsart}
\usepackage{latexsym, amsfonts, amsmath, amsthm, amssymb,amscd}
\usepackage[dvips]{graphicx}
\usepackage{eepic}
\usepackage{color}
\def\blue{\textcolor{blue}}
\def\stein{{Steingr\'{\i}msson}}
\def\os{\mathcal{OS}}
\def\tc{\mathcal{TC}}

\def\newterm#1{{\sl #1}\/}
\def\R{{\mathcal R}}

\def\O{{\mathcal O}}

\def\C{\mathcal{C}}
\def\S{\mathcal{S}}
\def\T{\mathcal{T}}
\def\N{\mathbb{N}}
\def\P{{\mathcal P}}
\def\OP{{\mathcal O}{\mathcal P}}

\def\D{{\mathcal D}}

\def\Inv{\mathop {\rm Inv}}
\def\inv{\mathop {\rm inv}}

\def\maj{\mathop {\rm maj}}
\def\comaj{\mathop {\rm comaj}}

\newcommand\eul{\operatorname{eul}}
\newcommand\mah{\operatorname{mah}}
\newcommand\Mah{\operatorname{Mah}}

\newcommand\rob{\operatorname{rob}}

\newcommand\cls{{\operatorname{cls}}}
\newcommand\opb{{\operatorname{opb}}}

\newcommand\sbg{{\operatorname{sb}}}
\newcommand\los{{\operatorname{los}}}
\newcommand\ros{{\operatorname{ros}}}
\newcommand\rcb{{\operatorname{rcb}}}
\newcommand\lcb{{\operatorname{lcb}}}
\newcommand\lob{{\operatorname{lob}}}
\newcommand\rcs{{\operatorname{rcs}}}
\newcommand\lcs{{\operatorname{lcs}}}

\newcommand\rsb{{\operatorname{rsb}}}
\newcommand\MAK{{\textsc{MAK}}}

\newcommand\des{\mathop{\rm des}}
\newcommand\mak{{\textsc{MAK}}}

\newcommand\Maj{\mathop{ \rm Maj}}
\newcommand\maf{\textsc{MAG}}
\newcommand\bDes{\mathop{ \rm bDes}}

\newcommand\stat{\mathop{ \textsc{STAT}}}

\newcommand\INV{\textsc{INV}}
\newcommand\MAJ{\textsc{MAJ}}

\newtheorem{thm}{Theorem}[section]
\newtheorem{prop}[thm]{Proposition}

\newtheorem{exam}{Example}[section]
\newtheorem{rmk}[thm]{Remark}
\newtheorem{Remark}[thm]{Remark}

\def\pf{\noindent {\it Proof.} }
\makeatletter
\newfont{\footsc}{cmcsc10 at 8truept}
\newfont{\footbf}{cmbx10 at 8truept}
\newfont{\footrm}{cmr10 at 10truept}
\makeatother 

\newtheorem{defn}[thm]{Definition}

\newtheorem{lemma}[thm]{Lemma}

\numberwithin{equation}{section}

\def\h{\operatorname{y}}

\def\open{\operatorname{open}}
\def\clos{\operatorname{clos}}

\def\los{\operatorname{los}}
\def\lcs{\operatorname{lcs}}
\def\ros{\operatorname{ros}}
\def\rcs{\operatorname{rcs}}
\def\lsb{\operatorname{lsb}}
\def\rsb{\operatorname{rsb}}
\def\lob{\operatorname{lob}}
\def\lcb{\operatorname{lcb}}
\def\rob{\operatorname{rob}}
\def\rcb{\operatorname{rcb}}
\def\bInv{\operatorname{bInv}}

\def\cbInv{\operatorname{cbInv}}
\def\bMaj{\operatorname{bMaj}}
\def\cbMaj{\operatorname{cbMaj}}
\def\cinvLSB{\operatorname{cinvLSB}}
\def\cmajLSB{\operatorname{cmajLSB}}

\def\sb{\operatorname{sb}}

\def\type{\lambda}

\renewcommand\tilde{\widetilde}

\begin{document}
\title{Euler-Mahonian Statistics
On \\Ordered Set
Partitions (II)
}

\author{Anisse Kasraoui}
\address{Universit\'e de Lyon; Universit\'e Lyon 1; Institut Camille
Jordan CNRS UMR 5208; 43, boulevard du 11 novembre 1918,
        F-69622, Villeurbanne Cedex}
\email{anisse@math.univ-lyon1.fr}

\author{Jiang Zeng}
\address{Universit\'e de Lyon; Universit\'e Lyon 1; Institut Camille
Jordan CNRS UMR 5208; 43, boulevard du 11 novembre 1918,
        F-69622, Villeurbanne Cedex}
\email{zeng@math.univ-lyon1.fr}

\begin{abstract}
We study  statistics on ordered set partitions whose generating
functions are related to $p,q$-Stirling numbers of the second
kind. The main purpose of this paper is to provide bijective
proofs of all the conjectures
 of \stein (Arxiv:math.CO/0605670). Our basic idea is to encode
ordered partitions by a kind of path diagrams and explore the rich
combinatorial properties of the latter structure. We also give  a
partition version of MacMahon's
 theorem on the equidistribution of the statistics
 inversion number and major index on words.
\end{abstract}

\maketitle

{\small \tableofcontents }


 \vskip 0.5cm \noindent{\bf MR Subject
Classifications}: Primary 05A18; Secondary 05A15, 05A30.\\


\section{Introduction}
The systematic study of statistics on permutations and
 words has its origins in the work of MacMahon~\cite{Ma}.
 In this paper,
 we will consider MacMahon's three statistics for a word $w$:
 the number of
\emph{descents} ($\des w$),  the number of \emph{inversions}
($\inv w$), and the \emph{major index} ($\maj w$). These are
defined as follows: A descent in a word $\pi=a_1a_2\cdots a_n$ is
an $i$ such that $a_i>a_{i+1}$,  an inversion is a pair $(i,j)$
such that $i<j$ and $a_i>a_j$, and the major index of $w$ is the
sum of the descents in $\pi$. The rearrangement class $R(w)$ of a
word $w=a_1a_2\cdots a_n$ is the set of
 all words obtained by permutating the letters of $w$.

Let ${\mathbf n}=(n_1,\ldots, n_k)$ be a sequence of non negative
integers and ${R}({\mathbf n})$  the rearrangement class  of the
word $1^{n_1}\ldots k^{n_k}$. Then MacMahon~\cite[Chap. 3]{And}
proved that
\begin{align}\label{eq:macmahon}
\sum_{w\in {R}({\mathbf n})}q^{\inv w} =\sum_{w\in {R}({\mathbf
n})}q^{\maj w} =\frac{(q;q)_{n_1+\cdots +n_k}}{(q;q)_{n_1}\cdots
(q;q)_{n_k}},
\end{align}
where $(x;q)_n=(1-x)(1-xq)\cdots (1-xq^{n-1})$. In particular, for
the symmetric group $S_n$ of $[n]:=\{1,2,\cdots,n\}$, we have
\begin{equation}\label{eq:MacMahon}
\sum_{\sigma\in\,\S_n}q^{\inv\,\sigma}=\sum_{\sigma\in\,\S_n}q^{\maj\,\sigma}=[n]!_q,
\end{equation}
where $[n]_q=1+q+\cdots q^{n-1}$ and $[n]_q!=[1]_q[2]_q\cdots
[n]_q$.

 Any statistic that
is equidistributed with $\des$ is said to be
\newterm{Eulerian}, while any statistic equidistributed with
$\inv$ is said to be \newterm{Mahonian}.  A bivariate statistic
$(\eul, \mah)$ is said to be a  \newterm{Euler-Mahonian} statistic
if  $\eul$ is Eulerian and $\mah$ is Mahonian.

An ordered  partition of $[n]$ is a sequence of disjoint and
nonempty subsets, called $\emph{blocks}$, whose union is $[n]$.
 The blocks of an ordered partition will be written as capital
letters separated by slashes, while elements of the blocks will be
set in lower case.
Thus an ordered partition of $[n]$ into $k$ blocks is written as
$\pi=B_1/B_2/\cdots /B_k$.
Clearly we can identify an
\emph{unordered} partition with an ordered partition
 by arranging the blocks in the increasing order of their minima, called
its  \emph{standard form}. For example, the partition $\pi$ of
$[9]$ consisting of the five blocks $ \{1,4,7\},\, \{2\},\,
\{3,9\},\, \{5\}\, \textrm{and}\, \{6,8\} $ will be written as
$\pi=1\;4\;7 /2/3\;9/5 /6\;8$. The set of all partitions of $[n]$
into $k$ blocks will be denoted by $\P_n^k$. It is well-known that
 the \emph{Stirling number of the second kind} $S(n,k)$ equals the cardinality  of $\P_n^k$.
Therefore,  if  $\OP_n^k$  denotes  the set of all ordered
partitions of $[n]$ into $k$ blocks, then $|\OP_n^k|=k!\,S(n,k)$.

For any ordered partition $\pi\in \OP_n^k$ there is a unique
partition $\pi_0=B_1/B_2/\cdots /B_k\in \P_n^k$ and  a unique
permutation $\sigma\in S_k$ such that
$\pi=B_{\sigma(1)}/B_{\sigma(2)}/\cdots /B_{\sigma(k)}$. In
parallel with notion of $\sigma$-restricted growth function in
\cite{Wa}, we shall call the corresponding partition $\pi$ a
\emph{$\sigma$-partition}. Let
 $\P_n^k(\sigma)$ be the set of all $\sigma$-partitions of
$[n]$ into $k$ blocks.
 For instance, $\pi= 6\;8 /5 /
1\;4\;7 / 3\;9/2\in\P_9^5(\sigma)$ with $\sigma=54132$.

Clearly,  for any $\sigma\in\S_k$ we have
$|\P_n^k(\sigma)|=|\P_n^k|=S(n,k)$ and
$\P_n^k=\P_n^k(\varepsilon)$ where $\varepsilon$ is the
\emph{identity permutation}.

 The $p,q$-Stirling numbers of
the second kind $S_{p,q}(n,k)$ were introduced in \cite{WW} by the
recursion:
\begin{align}\label{eq:stirling3}
 S_{p,q}(n,\, k)=
\left\{%
\begin{array}{ll}
    p^{k-1}S_{p,q}(n-1,\,k-1)+[k]_{p,q}\,S_{p,q}(n-1,\, k), & \hbox{if $0<k\leq n$;} \\
1, & \hbox{if $n=k=0$;} \\
0, & \hbox{otherwise.} \\\end{array}%
\right.
\end{align}
where
$$
[k]_{p,q}=p^{k-1}+p^{k-2}q+\cdots +pq^{k-1}+q^{k-1}.
$$

When $p$ or $q$ is set to 1, we obtain two usual $q$-Stirling
numbers of the second kind (see \cite{Gou}):
\begin{align}
S_q(n,k):=S_{q,1}(n,k)\quad \textrm{and} \quad \tilde S_{q}(n,\,
k):=S_{1,q}(n,k)=q^{-{k\choose 2}}S_q(n,k)
\end{align}

 Many authors (see e.g.
\cite{EhRe,GaRe,IKZ,KsZe,Mi,ReWa,Sa,Stein,Wa,WW,Wh}) have explored
the combinatorial aspects of these $q$-Stirling numbers. In
particular,  Wachs and White~\cite{WW} studied the combinatorial
interpretations of these $p,q$-Stirling numbers. In a sequel paper
Wachs \cite{Wa} extended some unordered partition interpretations
of $p,q$-Stirling number of the second kind
 to $\sigma$-partition statistics, although she used
 the \emph{restricted growth functions} instead of  set partitions.

A statistic $\stat$ on ordered set partitions is
\emph{Euler-Mahonian} if for any $n\geq k\geq 1$ its distribution over $\OP_n^k$ equals
$[k]_q!S_q(n,k)$, i.e.,
$$
\sum_{\pi\in\,\OP_n^k}q^{\stat\,\pi}=[k]_q!\,S_q(n,k).
$$

In \cite{Stein} Steingr\'{\i}msson conjectured several hard
Euler-Mahonian statistics on $\OP_n^k$. In a previous
paper~\cite{IKZ},  Ishikawa and the two current authors proved
half of the conjectures of Steingr\'{\i}msson~\cite{Stein} by
using the \emph{Matrix-transfer method} and determinant
computations.

The aim of this paper is to give a complete bijective approach to
Steingr\'{\i}msson's problem. In particular, we will not only
derive the results in \cite{IKZ} bijectively but also settle the
remaining half of the conjectures. In fact our bijective approach
yields also new results on $p,q$-Stirling numbers of Wachs and
White~\cite{WW} and the $\sigma$-partitions of Wachs~\cite{Wa}. As
we will show, one of our results generalizes MacMahon's
equidisribution result of inversion number and major index on
words.

Throughout this paper, we shall denote by $\mathbb{P}$ (resp.
$\N$, $\mathbb{Z}$) the set of positive integers (resp. non
negative integers, integers) and assume that $n$ and $k$ are two
fixed integers satisfying $n\geq k\geq1$. Furthermore, for any
integers $i_1, i_2,\ldots,i_k$, we denote by
$\{i_1,i_2,\ldots,i_k\}_<$ (resp. $\{i_1,i_2,\ldots,i_k\}_>$) the
increasing (resp. decreasing) arrangement of these integers.
\section{Definitions}
Let $B$ be a  finite subset of $\N$. The \emph{opener} of
$B$  is its least element while the
\emph{closer} of $B$  is its greatest
element. For $\pi\in\OP_n^k$, we will denote by $\open(\pi)$ and
$\clos(\pi)$ the sets of openers and closers of the blocks of $\pi$, respectively.
 The letters (integers) in $\pi$ are further divided into four classes:
\begin{itemize}
 \item \emph{singletons}: elements of the singleton blocks;
 \item \emph{strict openers}:  smallest elements  of the non singleton
 blocks;
 \item \emph {strict closers}: largest elements of the non singleton
 blocks;
\item \emph{transients}: all other elements, i.e., non extremal
elements of non singleton blocks.
\end{itemize}
The sets of strict openers, strict closers, singletons and
transients of $\pi$ will be denoted, respectively, by $\O(\pi)$,
$\C(\pi)$, $\S(\pi)$ and $\T(\pi)$. Obviously we have
$$\open(\pi)=\O(\pi)\cup\S(\pi),\quad
 \clos(\pi)=\C(\pi)\cup\S(\pi),\quad \S(\pi)=\open(\pi)\cap \clos(\pi).
 $$
The $4$-tuple $(\O(\pi),\C(\pi),\S(\pi),\T(\pi))$, denoted by
$\type(\pi)$, is called the \emph{type} of $\pi$.
 For example, if $\pi=3\,5/2\,4\,6/1/7\,8$, then
$\open(\pi)=\{1,2,3,7\}$, $\clos(\pi)=\{1,5,6,8\}$ and
$$
  \type(\pi)=(\{2,3,7\},\{5,6,8\},\{1\},\{4\}).
$$
Let $\OP_n^k(\lambda)$ be the set of ordered partitions in
$\OP_n^k$ of type $\lambda$.

Following Steingr\'{\i}msson \cite{Stein}, we now define a system of ten
inversion-like statistics on ordered set partitions. Most of them
have been studied in the case of set partitions. Note that the two
last statistics were essentially defined by Foata and Zeilberger
\cite{FoZe} for permutations. The reader is refered to
\cite{KsZe,Stein,Wa} for further informations of these statistics.

Given a partition $\pi=B_1/B_2/\cdots /B_k\in\OP_n^k$, let $w_i$ be
the index of the block (counting from left to right) containing $i$,
namely the integer $j$ such that $i\in B_j$. Then define ten
coordinate statistics as follows. For $1\leq i\leq n$, we let:
\begin{eqnarray*}
\los_i\,\pi &=& \# \{j \in \open\,\pi\, |\,j<i, \,w_j<w_i\},\\
\ros_i\,\pi &=& \# \{j \in \open\,\pi\, |\,j<i, \,w_j>w_i\}, \\
\lob_i\,\pi &=& \# \{j \in \open\,\pi\, |\,j>i, \,w_j<w_i\}, \\
\rob_i\,\pi &=& \# \{j \in \open\,\pi\, |\,j>i, \,w_j>w_i\}, \\
\lcs_i\,\pi &=& \# \{j \in \clos\,\pi\, |\,j<i, \,w_j<w_i\}, \\
\rcs_i\,\pi &=& \# \{j \in \clos\,\pi\, |\,j<i, \,w_j>w_i\},\\
\lcb_i\,\pi &=& \# \{j \in \clos\,\pi\, |\,j>i,\,w_j<w_i\},\\
\rcb_i\,\pi &=& \# \{j \in \clos\,\pi\, |\,j>i,\,w_j>w_i\}.
\end{eqnarray*}
Moreover, let $\rsb_i\,\pi$ (resp. $\lsb_i\,\pi$) be the number of
blocks \textrm{B} in $\pi$ to the right (resp. left) of the block
containing $i$ such that the opener of \textrm{B} is smaller than
$i$ and the closer of \textrm{B} is greater than $i$. Remark that
$\lsb_i$ and $\rsb_i$ are each equal to the difference of two of the
first eight statistics. Namely, it is easy to see that
\begin{equation}\label{eq:decompo lsb,rsb}
\lsb_i=\los_i-\lcs_i=\lcb_i-\lob_i\quad\text{and}\quad
\rsb_i=\ros_i-\rcs_i=\rcb_i-\rob_i.
\end{equation}
Then define the statistics $\ros$, $\rob$, $\rcs$, $\rcb$, $\lob$,
$\los$, $\lcs$, $\lcb$, $\lsb$ and $\rsb$ as the sum of their
coordinate statistics, e.g.
$$
\ros=\ros_{1}+\cdots \ros_{n}.
$$
For any partition $\pi$ we can define the restrictions of these statistics
on openers and non openers
$\ros_{\os}$, $\ros_{\tc}$, ..., $\rsb_{\os}$ and $\rsb_{\tc}$, e.g.
\begin{align}\label{def:partialstat}
{\ros}_{\os}\,\pi=\sum_{i\in\,
\O\cup\S(\pi)}{\ros}_i\,\pi\quad\text{and}\quad
{\ros}_{\tc}\,\pi=\sum_{i\in\, \T\cup \C(\pi)}{\ros}_i\,\pi.
\end{align}

\begin{rmk}
Note that $\ros$ is the abbreviation of "right, opener, smaller",
while $\lcb$ is the abbreviation of "left, closer, bigger", etc.
\end{rmk}

 As an example, we give here the values of the coordinate  statistics computed
on $\pi=  6\;8 /5  / 1\;4\;7  / 3\;9/2$:
$$
\begin{array}{cccccccccc}
\pi=&  6\;8 & / & 5 & / & 1\;4\;7 & / & 3\;9 & / & 2\\
&&&&&&&&&\\
   \los_i: & 0\;0 & / & 0 & / & 0\;0\;2 & / & 1\;3 & / & 1 \\
   \ros_i: & 4\;4 & / & 3 & / & 0\;2\;2 & / & 1\;1 &  /& 0 \\
   \lob_i: & 0\;0 & / & 1 & / & 2\;2\;0 & / & 2\;0 & / & 3 \\
   \rob_i: & 0\;0 & / & 0 & / & 2\;0\;0 & / & 0\;0 & /& 0 \\
   \lcs_i: & 0\;0 & / & 0 & / & 0\;0\;1 & / & 0\;3 & / & 0 \\
   \rcs_i: & 2\;3 & / & 1 & / & 0\;1\;1 & / & 1\;1 & / & 0 \\
   \lcb_i: & 0\;0 & / & 1 & / & 2\;2\;1 & / & 3\;0 & / & 4 \\
   \rcb_i: & 2\;1 & / & 2 & / & 2\;1\;1 & / & 0\;0 & / & 0 \\
   \lsb_i: & 0\;0 & / & 0 & / & 0\;0\;1 & / & 1\;0 & / & 1 \\
   \rsb_i: & 2\;1 & / & 2 & / & 0\;1\;1 & / & 0\;0 & / & 0 \\
\end{array}
$$
It follows
that ${\ros}_{\os}\,\pi=8$ and ${\rsb}_{\tc}\,\pi=3$.

The following result is due to Wachs and White~\cite[Cor. 5.3]{WW}.
\begin{equation}\label{eq:WW}
\sum_{\pi\in\P_n^k}p^{\rcb\,\pi}q^{\lsb\,\pi}=S_{p,q}(n,k).
\end{equation}

Let $\pi=B_1/B_2/\cdots /B_k \in \OP_n^k$. Define a \emph{partial
order} $\succ$ on blocks $B_i$'s as follows~: $B_i\succ B_j$ if
all the letters of $B_i$ are greater than those of $B_j$, i.e., if
$\min(B_i)>\max(B_j)$.
 A  \emph{block inversion} in $\pi$ is a pair $(i,j)$ such that $i<j$ and
 $B_{i} \succ B_{j}$. We denote by  $\bInv\, \pi$ the number of block
inversions in $\pi$.  A \emph{block descent} is an integer $i$
such that $B_i \succ B_{i+1}$.
 The \emph{block major index} of $\pi$, denoted by $\bMaj\,\pi$, is the sum of the block
descents in $\pi$. We also define their complementary
counterparts:
\begin{align}
\cbInv = {k\choose2}-\bInv\,\qquad\cbMaj = {k\choose2}-\bMaj,
\end{align}
and two composed statistics:
\begin{align}\label{def:cinv-cmaj}
\begin{gathered}
\cinvLSB=\lsb + \cbInv + {k\choose2},\\
\cmajLSB=\lsb + \cbMaj + {k\choose2}.
\end{gathered}
\end{align}

For any $\sigma$-partition $\pi\in\P_n^k(\sigma)$ with $\sigma\in\S_k$, we
can  define two natural statistics~:
\begin{equation}\label{eq:inv,maj-perm-part}
\Inv\,\pi=\inv\,\sigma\quad\text{and}\quad \Maj\,\pi=\maj\,\sigma.
\end{equation}

The following two statistics are the (ordered)
partition analogues of their counterparts in permutations \cite[p.13]{Stein}:

\begin{align}\label{def:mak-mak'}
\begin{gathered}
\MAK= \ros +\lcs,\qquad
{\MAK}'= \lob +\rcb.
\end{gathered}
\end{align}
In what follows we will also denote by $\maf$ any of these two statistics, i.e.,
\begin{align}\label{def:mag}
\maf\in \{\MAK, \MAK'\}.
\end{align}


\section{Main results}
We first present  our main result on the equidistribution of some
inversion like statistics  on $\sigma$-partitions.

\begin{thm}\label{thm:sym}
 For any $\sigma\in\S_k$, the triple statistics
$(\MAK+\bInv,\MAK'+\bInv,\cinvLSB)$ and
$(\MAK'+\bInv,\MAK+\bInv,\cinvLSB)$
 are equidistributed on
$\P_n^k(\sigma)$.
Moreover,
\begin{align}\label{eq:KZWWbis1General}
\sum_{\pi\in
\P_n^k(\sigma)}p^{(\maf+\bInv)\, \pi}q^{\cinvLSB\,\pi}
=q^{k(k-1)}\left(\frac{p}{q}\right)^{\inv\sigma}S_{p,q}(n,k).
\end{align}
\end{thm}
Note that \eqref{eq:KZWWbis1General} gives  two
$\sigma$-extensions of \eqref{eq:WW}, which is the
$\sigma=\epsilon$ case of \eqref{eq:KZWWbis1General}.

As $\inv$ is a Mahonian statistic on $S_k$,   summing the two
sides over all the permutations $\sigma$ in $S_k$, we derive
immediately  from \eqref{eq:KZWWbis1General} the main result of
\cite{IKZ}, of which the second part   was conjectured by
Steingr\'{\i}msson~\cite{Stein}.

\begin{thm}[Ishikawa et al.]\label{thm:IKZ}
We have
\begin{align}\label{eq:Mah1}
\sum_{\pi\in
\OP_n^k}p^{(\maf+\bInv)\pi}q^{\cinvLSB \pi}=q^{k\choose
2}[k]_{p,q}!\,S_{p,q}(n,k).
\end{align}
In particular the three statistics
\begin{align}
\mak +\bInv, \quad \mak'+\bInv,\quad\cinvLSB\label{invstat}
\end{align}
are Euler-Mahonian on $\OP_n^k$.
\end{thm}

For the  major like statistics we have the following
equidistribution result on ordered partitions with a fixed type.
\begin{thm}\label{thm:refinesteinconj}
 The triple statistics
$(\mak+\bMaj, \mak'+\bMaj, \cmajLSB)$ and
$(\mak+\bInv, \mak'+\bInv, \cinvLSB)$
are equidistributed on $\OP_n^k(\lambda)$ for  any partition type $\lambda$.
\end{thm}

Combining Theorems~\ref{thm:IKZ} and \ref{thm:refinesteinconj} we
derive immediately  the following result, of which   the first
part  was  conjectured by Ishikawa et al. \cite[Conjecture 6.2]{IKZ} while the
second part was originally conjectured by
Steingr\'{\i}msson~\cite{Stein}.

\begin{thm}\label{thm:cor1}
We have
\begin{align}\label{eq:Mah1}
\sum_{\pi\in
\OP_n^k}p^{(\maf+\bMaj)\pi}q^{\cmajLSB \pi}=q^{k\choose
2}[k]_{p,q}!\,S_{p,q}(n,k).
\end{align}
In particular  the three statistics
 \begin{align}
\mak+\bMaj, \quad \mak'+\bMaj,\quad  \cmajLSB\label{majstat}
\end{align}
are Euler-Mahonian on $\OP_n^k$.
\end{thm}

For any set partition $\pi\in \P_n^k$, denote by $\R(\pi)$ the
rearrangement class of $\pi$, i.e., if $\pi=B_1/B_2/\cdots
 /B_k$, then
 $$
 \R(\pi)=\{B_{\sigma(1)}/B_{\sigma(2)}/\cdots /B_{\sigma(k)}\;|\;\sigma\in\S_k \}.
 $$
 For instance, if
$\pi=1\;4/2\;3/5$, then
$$
\R(\pi)=\{1\;4/2\;3/5 \,,\, 1\;4/5/2\;3 \,,\, 5/2\;3/1\;4 \,,\,
2\;3/1\;4/5 \,,\,\, 2\;3/5/1\;4 \,,\, 5/1\;4/2\;3\}.
$$
It is clear that  $|\R(\pi)|=k!$ for any  partition $\pi$ with $k$
blocks.

Introduce first two analogues of inversion numbers and major index
on $\OP_n^k$:
\begin{align}\label{eq:inv,maj}
\begin{gathered}
\INV={\rsb}_{\os}+\bInv,\\
\MAJ={\rsb}_{\os}+\bMaj.
\end{gathered}
\end{align}
The statistic $\INV$ is just a rewording of the inversion number
$\inv$ in \eqref{eq:inv,maj-perm-part}. Indeed, it is easy to see
that $\bInv={\rcs}_{\os}$ and $\Inv={\ros}_{\os}$, therefore
\begin{equation}\label{eq:defnalternative de INV}
\INV={\rsb}_{\os}+{\rcs}_{\os}=\ros_{\os}=\Inv.
\end{equation}

Our last result is  a non trivial extension of MacMahon's identity
\eqref{eq:macmahon} to rearrangement class of an arbitrary
partition.
\begin{thm}\label{thm:extMac}
For any $\pi\in\P_n^k$, the statistics $\INV$ and $\MAJ$ are
equidistributed on $\R(\pi)$ and
\begin{equation}\label{extmac}
\sum_{\pi\in\R(\pi)}q^{\MAJ\,
\pi}=\sum_{\pi\in\R(\pi)}q^{\INV\,\pi}=[k]_q!.
\end{equation}

\end{thm}
 To show that  \eqref{extmac} implies
  MacMahon's formula \eqref{eq:macmahon}, we consider the rearrangement class of
  a special set partition
  as follows.
Let  $N_i=n_1+\cdots +n_i$ for $i=1,\ldots, k$ and
$$
\Pi=\pi_{11}\ldots \pi_{1n_1}\pi_{21}\ldots \pi_{2n_2}\ldots
\pi_{k1}\ldots \pi_{kn_k}
$$ be the partition of $[2N_k]$ consisting
of the doubletons:
\begin{align*}
\pi_{ij}:=\{2N_{i-1}+j,\; 2N_{i-1}+n_i+j\},\qquad (1\leq i\leq
k,\; 1\leq j\leq n_i,\; N_0=0).
\end{align*}

 It is readily  seen that each
$\pi\in \R(\Pi)$ can be identified with a pair $(w, (\pi_1,\ldots,
\pi_k))$ where $w\in {R}({\mathbf n})$ is the word obtained from
$\pi$ by substituting each $\pi_{ij}$ by $i$ for $1\leq i\leq k$
and $\pi_i$ is the word obtained from $\pi$ by deleting all the
$\pi_{lj}$ for $l\neq i$. For example, if $k=3$,
 $n_1=3$, $n_2=2$ and $n_3=3$ then
\begin{align*}
\Pi&=\{1,4\}\,\{2,5\}\,\{3,6\}
\,\{7,9\}\,\{8,10\}\,\{11,14\}\,\{12,15\}\,\{13,16\}.
\end{align*}
Let $\pi\in \R(\Pi)$ be the ordered partition:
$$\pi=\{7,9\}\,\{2,5\}\,\{11,14\}\,\{1,4\}\,\{13,16\}\,\{12,15\}\,\{3,6\}
\,\{8,10\}.
$$
Then $\pi\mapsto (w, (\pi_1,\pi_2,\pi_3))$ with
$w=2\,1\,3\,1\,3\,3\,2\,1$ and
$$
\pi_1=\{2,5\}\,\{1,4\}\,\{3,6\},\;\pi_2=\{7,9\}\,\{8,10\},\;
\pi_3=\{11,14\}\,\{13,16\}\,\{12,15\}.
$$

Note that for any $i,j,l,r$ such  that $1\leq i,j\leq k$, $1\leq
l\leq n_i$ and $1\leq r\leq n_j$,
$$
\pi_{il}\prec\pi_{jr} \Longleftrightarrow i<j,
$$
therefore
 \begin{align*}
\bMaj\pi=\maj w,\qquad \bInv\pi=\inv w
 \end{align*}
and
 $$\rsb_{\O\S}\pi=\rsb_{\O\S}\pi_1+\cdots
+\rsb_{\O\S}\pi_k.
$$
It then follows from \eqref{eq:inv,maj} that
 \begin{align}
\sum_{\pi\in \R(\Pi)}q^{\MAJ\pi}&= \left(\sum_{w\in {\R}({\mathbf
n})}q^{\maj w}\right)
\prod_{i=1}^k\sum_{\pi\in \R(\pi_i)}q^{\rsb_{\O\S}\pi},\label{mac1}\\
\sum_{\pi\in \R(\Pi)}q^{\INV\pi}&= \left(\sum_{w\in {\R}({\mathbf
n})}q^{\inv w}\right) \prod_{i=1}^k\sum_{\pi\in
\R(\pi_i)}q^{\rsb_{\O\S}\pi} .\label{mac2}
 \end{align}
As $\rcs_\os(\pi)=0$ for any $\pi\in \R(\pi_i)$ we have
$\INV={\rsb}_{\os}$ on $\R(\pi_i)$. Hence, by
Theorem~\ref{thm:extMac},
 \begin{align*}\label{eq:maj}
\sum_{\pi\in \R(\pi_i)}q^{\rsb_{\O\S}\pi}=[n_i]_q!,\qquad 1\leq
i\leq k,
 \end{align*}
 and
 $$
\sum_{\pi\in \R(\Pi)}q^{\INV\pi}=[N_k]_q!.
 $$
 MacMahon's formula \eqref{eq:macmahon} follows then from \eqref{mac1} and \eqref{mac2} by
invoking Theorem~\ref{thm:extMac}.

 The  rest of this paper is organized  as follows.
In Section~4 we shall  present the first path diagrams encoding of
ordered partitions $\Phi: \Delta_n^k\longrightarrow  \OP_n^k$,
which was introduced in \cite{IKZ}. We will construct an
involution on path diagrams $\varphi: \Delta_n^k\longrightarrow
\Delta_n^k$ in Section~5 and a bijection $\Gamma_\sigma$ from set
partitions to ordered set partitions in Section~6. We prove
Theorem~3.1 in section~7. To deal with major like statistics a
second path diagram encoding $\Psi:\Delta_n^k\longrightarrow
\OP_n^k$ will be given in Section~8. Then, in Section~9 we prove
Theorem~3.4 by using the mapping $\Upsilon:=\Psi \circ
\Phi^{-1}:\OP_n^k\to \OP_n^k$. Finally we prove Theorem~3.6 in
Section~3.6 and conclude the paper with some further remarks.

\section{The first path diagram encoding $\Phi$ of  ordered partitions}
As shown in \cite{IKZ}, we can define
the notion of \emph{trace} or \emph{skeleton} for ordered partitions. To
this end, adjoin to $\mathbb{P}$ a symbol $\infty$ such that $i<\infty$ for any positive
integer $i$. The restriction of a subset $B$ of $\mathbb{P}$ on
$[i]$, namely $B\cap [i]$, is said to be
\begin{itemize}
\item \emph{empty} if $i<\min B$,
\item \emph{active} if $\min B\leq i<\max B$,
\item \emph{complete} if $\max B\leq i$.
\end{itemize}

Let $\pi\in\OP_n^k$.  For $1\leq i\leq n$ the \emph{$i$-trace} or
$i$-skeleton $T_i$  of $\pi$ is obtained by restricting each block
on $[i]$ and deleting  empty blocks. By convention, we  add a
symbol $\infty$ at the end of  each active block. By convention, the latter
 is still called block.

Clearly  one can characterize a partition and the statistics $\ros$ and
$rsb$  by using its traces. More precisely, for any $1\leq i\leq n$,
$\ros_i\pi$
(resp. $\rsb_i\pi$)  equals the number of blocks (resp. active blocks) to the right of
the block containing $i$ in the $i$-trace $T_i$ of $\pi$.
\begin{exam}
If $\pi=3\;5\;7/1\;4\;10/9/6/2\;8$, then the  $7$-trace of $\pi$ is
$$
T_7=3\;5\;7/1\;4\,\infty/6/2\,\infty,
$$
with  $\ros_7\pi=3$ and $\rsb_7\pi=2$.
\end{exam}

 A useful
way to describe the traces of a partition is to draw a \emph{path
diagram}.

\begin{defn} A path of depth $k$ and length $n$ is a sequence
$w=(w_0, w_1, \ldots, w_n)$ of points in $\mathbb{N}^2$ such that
$w_0=(0,0)$, $w_n=(k,0)$ and the $i$-th ($1\leq i\leq n$) step
 $(w_{i-1}, w_i)$ must be one
of the following four types:
\begin{itemize}
\item North, i.e.,  $w_i-w_{i-1}=(0,1)$,
\item East, i.e., $w_i-w_{i-1}=(1,0)$,
\item South-East, i.e., $w_i-w_{i-1}=(1,-1)$,
\item Null, i.e., $w_i-w_{i-1}=(0,0)$ and $y_i>0$ if $w_i=(x_i,y_i)$.
\end{itemize}
The \emph{abscissa} and \emph{ordinate} of $w_{i-1}$ are called
the \emph{abscissa} and  \emph{height} of the $i$-th step of $w$
and denoted  by $x_i(w)$ and $\h_i(w)$, respectively. The set of
all paths of depth $k$ and length $n$ will be denoted by
$\Omega_n^k$.
\end{defn}

We can visualize a path $w$ by drawing a \emph{segment} or
\emph{loop} from $w_{i-1}$ to $w_i$ in the $xy$ plane. For
instance, the path
$$w=((0,0),\,(0,1),\,(0,2),\,(0,3),\,(0,3),\,
(0,3),\,(1,3),\,(2,2),\,(3,1),\,(4,1),\,(5,0)),
$$
 is illustrated  in  Figure~1.
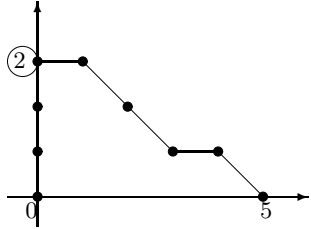
\begin{figure}[h]\
\begin{center}
{\setlength{\unitlength}{0.4mm}
\begin{picture}(90,60)(0,0)
\put(0,0){\circle*{3}}\put(0,15){\circle*{3}}
\put(0,30){\circle*{3}}\put(0,45){\circle*{3}}
\put(15,45){\circle*{3}}
\put(30,30){\circle*{3}}\put(45,15){\circle*{3}}
\put(60,15){\circle*{3}}\put(75,0){\circle*{3}}
\put(-10,0){\vector(1,0){100}}\put(0,-10){\vector(0,1){75}}
\put(0,0){\linethickness{0.2mm}\line(0,1){45}}
\put(0,45){\linethickness{0.2mm}\line(1,0){15}}
\put(15,45){\linethickness{0.2mm}\line(1,-1){30}}
\put(45,15){\linethickness{0.2mm}\line(1,0){15}}
\put(60,15){\linethickness{0.2mm}\line(1,-1){15}}
\put(-5,45){\linethickness{0.2mm}\circle{10}}\put(-8,43){\scriptsize2}
\put(-4,-7){\scriptsize 0}\put(74,-7){\scriptsize5}
\end{picture}}
\end{center}
\caption{A path in $\Omega_{10}^5$ with two successive Null steps
from $(0,3)$ to $(0,3)$. }
\end{figure}

For the reason which will be clear later, the sets of indices of North, South-East,
East and Null steps of a path $w$ will be denoted by $\O(w)$,
$\C(w)$, $\S(w)$ and $\T(w)$, respectively. The $4$-tuple
$(\O(w),\C(w),\S(w),\T(w))$, denoted by $\type(w)$, is called the
\emph{type} of $w$.
  For instance, if $w$ is the path represented in Figure~1, then
$$
  \type(w)=(\{1,2,3\},\{7,8,10\},\{6,9\},\{4,5\}).
$$
\begin{defn}
A \emph{path diagram of depth $k$ and length $n$} is a pair
$(w,\gamma)$, where $w$ is a path in $\Omega_n^k$ and
$\gamma=(\gamma_i)_{1\leq i\leq n}$ is a sequence of integers such
that
\begin{itemize}
\item $0\leq \gamma_i\leq y_i(w)-1$ if the $i$-th step of $w$ is
Null or South-East,
\item  $0\leq \gamma_i\leq x_i(w)+y_i(w)$ if
the $i$-th step of $w$ is North or East.
\end{itemize}
Let $\Delta_n^k$ be the set of path diagrams of depth $k$
and length $n$.
\end{defn}

For any statement $A$ we denote by $\chi(A)$ the character
function of $A$, that means $\chi(A)=1$ if $A$ is true and
$\chi(A)=1$  if $A$ is false.
\begin{lemma} \label{usefullemma}
Let $w\in \Omega_n^k$ and $(w,\gamma)\in \Delta_n^k$ with
$\O\cup\S(w)=\{i_1, i_2,\ldots,i_k\}_<$. Then
\begin{itemize}
\item[(i)] $x_{i}+y_{i}=x_{i-1}+y_{i-1}+\chi(i\in \O\cup\S(w))$ for $i=1,\ldots,n$.
\item[(ii)] $x_{i_j}+y_{i_j}+1=j$  for $ j=1,2,\ldots, k$.
\end{itemize}
 \end{lemma}
 \begin{proof}
Since $x_1(w)+y_1(w)=0$ and each step is one of the four kinds:
 $(0,1)$, $(1,0)$, $(0,0)$ and $(1,-1)$,
 the sum $x_i+y_i$ increases by one if and only if the
$i$-th step is North or South. This yields (i), while  (ii) is a
direct consequence of (i).

\end{proof}

One can encode ordered partitions by
path diagrams. The following  important bijection
$\Phi: \Delta_n^k \to \OP_n^k$
was introduced in \cite{IKZ}.

\noindent{\bf Algorithm $\Phi$}: Starting from  a path diagram
$h=(w,\gamma)$ in $\Delta_n^k$, we obtain $\Phi(h)=\pi$ by
constructing recursively all the
 $i$-traces $T_i$ ($1\leq i\leq n$) of $\pi$, i.e., such that $\pi=T_n$.
 By convention  $T_0=\emptyset$.
 Assume that we have constructed
$T_{i-1}=B_{i_1}/B_{i_2}/\cdots /B_{i_{\ell}}$  such that  $T_{i-1}$
has $\h_i(w)$ active blocks and $x_{i}(w)$ complete blocks. Label
the slots before $B_{i_1}$, between $B_{i_j}$ and $B_{i_{j+1}}$,
for $1\leq j\leq \ell-1$, and after $B_{i_{\ell}}$ from left to right by
$\ell,\ldots,1, 0$, while the active blocks of $T_{i-1}$ are labeled
from  left to right by $\h_i(w)-1,\ldots, 1,0$. Extend  $T_{i-1}$
to $T_{i}$ as follows:
 \begin{itemize}
\item If the $i$-th step of $w$ is North (resp. East),
 then create an active  block (resp. singleton) with $i$ at the  slot with label $\gamma_i$;
\item If the $i$-th step of $w$ is Null (resp. South-East),
 then insert $i$ (resp. replace $\infty$
by $i$) in the  active block with label $\gamma_i$.
\end{itemize}

Since $x_{n+1}(w)=k$ and $\h_{n+1}(w)=0$,
the $n$-trace  $T_n$ is a partition in  $\OP_n^k$.

\begin{figure}\label{fig:Psi}
\begin{center}
{\setlength{\unitlength}{0.4mm}
\begin{picture}(90,330)(0,-280)

\put(-30,30){\makebox(3,4){$h\;=$}}
\put(0,0){\circle*{3}}\put(0,15){\circle*{3}}\put(0,30){\circle*{3}}\put(0,45){\circle*{3}}\put(15,45){\circle*{3}}
\put(30,30){\circle*{3}}\put(45,15){\circle*{3}}\put(60,15){\circle*{3}}\put(75,0){\circle*{3}}
\put(-10,0){\vector(1,0){100}}\put(0,-10){\vector(0,1){75}}\put(0,0){\linethickness{0.2mm}\line(0,1){45}}\put(0,45){\linethickness{0.2mm}\line(1,0){15}}
\put(15,45){\linethickness{0.2mm}\line(1,-1){30}}\put(45,15){\linethickness{0.2mm}\line(1,0){15}}\put(60,15){\linethickness{0.2mm}\line(1,-1){15}}
\put(-5,45){\linethickness{0.2mm}\circle{10}}

\put(-4,7){\scriptsize0}\put(-4,22){\scriptsize0}\put(-4,34){\scriptsize2}
\put(-20,46){\scriptsize1,2}\put(7,47){\scriptsize3}\put(24,40){\scriptsize2}
\put(39,26){\scriptsize0}\put(52,17){\scriptsize1}\put(69,11){\scriptsize0}

\put(35,-5){\vector(0,-1){35}}\put(25,-29){\textbf{$\Phi$}}

\put(0,-250){\makebox(100,160){$ \scriptsize{
\begin{array}{cccccl}
i & \quad\textrm{step}_i\quad& \gamma_i &&& \quad T_i\\
\hline\\
0 &&  &&& \emptyset \\
\hline\\
1 & North& 0 &&& \textbf{1}\,\infty \\
\hline\\
2 & North& 0 &&& 1\,\infty/\textbf{2}\,\infty \\
\hline\\
3 & North& 2 &&& \textbf{3}\,\infty/1\,\infty /2\,\infty \\
\hline\\
4 & Null & 1 &&& 3\,\infty/1\;\textbf{4}\,\infty/2\,\infty \\
\hline\\
5 &Null & 2 &&& 3\;\textbf{5},\infty/1\;4,\infty/2\,\infty \\
\hline\\
6 &East & 3 &&& \textbf{6}/3\;5\,\infty/1\;4\,\infty/2\,\infty \\
\hline\\
7 & South\text{-}East& 2 &&& 6/3\;5\;\textbf{7}/1\;4\,\infty/2\,\infty \\
\hline\\
8 & South\text{-}East& 0 &&& 6/3\;5\;7/1\;4\,\infty/2\;\textbf{8} \\
\hline\\
9 & East & 1 &&& 6/3\;5\;7/1\;4\;\infty/\textbf{9}/2\;8\\
\hline\\
10 & South\text{-}East & 0
&&& 6/3\;5\;7/1\;4\;\textbf{10}/9/2\;8.\\
\hline\\
\end{array}}
$}}
\end{picture}}
\end{center}
\caption{the step by step construction of $\Phi(h)$}
\end{figure}

\begin{exam}\label{exam:Phi}
 Consider the path diagram $h=(w,\gamma)\in\Delta_{10}^5$, where
 $w$
is the path in Figure~1 and $\gamma=(0,0,2,1,2,3,2,0,1,0)$,
 then $\Phi(h)=6/3\;5\;7/1\;4\;10/9/2\;8$. The step
by step construction of $\Phi(h)$ is given in Figure~2, where the
$i$-th step is labeled with  $\gamma_i$ for $1\leq i\leq 10$.
\end{exam}

The main properties of $\Phi$ are given in the following theorem of Ishikawa et al~\cite{IKZ}.
\begin{thm}\label{thmPhi}
 The bijection $\Phi: (w,\gamma)\mapsto \pi$ has the following properties:
$\type(w)=\type(\pi)$ and for $1\leq i\leq n$,
 $$
\gamma_i=\left\{%
\begin{array}{ll}
    \ros_i(\pi), & \hbox{if $i\in \O(\pi)\cup\S(\pi)$;} \\
   \rsb_i(\pi), & \hbox{if $i\in \T(\pi)\cup\C(\pi)$.} \\
\end{array}%
\right.
$$
\end{thm}
It follows that
\begin{equation}\label{eq:thmPsi}
{\ros}_{\os}\,\pi=\sum_{i\in\O\cup\S(w)}\gamma_i \quad\text{and}\quad
{\rsb}_{\tc}\,\pi=\sum_{i\in\T\cup\C(w)}\gamma_i.
\end{equation}

\section{Involution $\varphi$ on path diagrams}
For  any path $w\in\Omega_n^k$,
one can define a natural bijection between the North steps and
South-East steps  of $w$.
 More precisely, let
$$\O(w)=\{o_1,o_2,\cdots,o_r\}_{<}\quad \textrm{and}\quad
\C(w)=\{c_1,c_2,\cdots,c_r\}_{<}.
$$
 We define the
permutation $\sigma\in S_r$, called the \emph{associated permutation}
of $w$,
 as follows:
 Suppose the height of the $i$-th North
step of $w$ (i.e. the $o_i$-th step of w) is $t$. Since
$w_0=(0,0)$ and $w_n=(k,0)$, there must exist a South-East step of
height $t+1$ to the right of the $o_i$-th step. If the first such
step is the $j$-th South-East step (i.e. the  $c_j$-th
step of $w$), set $\sigma(i)=j$. Since there is at least one
South-East step of height $t+1$ between any two North steps of
height $t$ the mapping $\sigma$ is an injection and then  a
bijection for $\O(w)$ and $\C(w)$ have the same cardinality.

For instance, if $w$ is the path represented below,
 we have $\sigma(1)=4$, $\sigma(2)=2$, $\sigma(3)=1$ and
$\sigma(4)=3$. The construction of $\sigma$ is illustrated in
Figure~3.

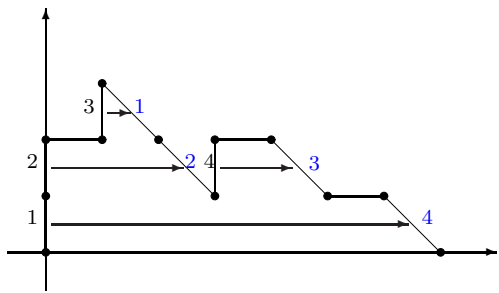
\begin{figure}[h]\label{fig:associated permutation}
\begin{center}
{\setlength{\unitlength}{0.5mm}
\begin{picture}(90,60)(0,0)
\put(0,0){\circle*{2}}\put(0,15){\circle*{2}}\put(0,30){\circle*{2}}\put(15,30){\circle*{2}}\put(15,45){\circle*{2}}
\put(30,30){\circle*{2}}\put(45,15){\circle*{2}}\put(45,30){\circle*{2}}\put(60,30){\circle*{2}}\put(75,15){\circle*{2}}
\put(90,15){\circle*{2}}\put(105,0){\circle*{2}}

\put(-10,0){\vector(1,0){130}}\put(0,-10){\vector(0,1){75}}
\put(0,0){\linethickness{0.2mm}\line(0,1){30}}\put(0,30){\linethickness{0.2mm}\line(1,0){15}}
\put(15,30){\linethickness{0.2mm}\line(0,1){15}}\put(15,45){\linethickness{0.2mm}\line(1,-1){30}}
\put(45,15){\linethickness{0.2mm}\line(0,1){15}}\put(45,30){\linethickness{0.2mm}\line(1,0){15}}
\put(60,30){\linethickness{0.2mm}\line(1,-1){15}}\put(75,15){\linethickness{0.2mm}\line(1,0){15}}
\put(90,15){\linethickness{0.2mm}\line(1,-1){15}}

\put(-5,7.5){\tiny 1}\put(-5,22.5){\tiny{2}}
\put(10,37){\tiny3}\put(42,22.5){\tiny 4}

\put(23.5,37){\blue{\tiny 1}}\put(37,22.5){\blue{\tiny{2}}}
\put(70,22){\blue{\tiny 3}}\put(100,7.5){\blue{\tiny4}}

\put(1.5,7.5){\vector(1,0){95}}\put(1.5,22.5){\vector(1,0){35}}
\put(16.5,37){\vector(1,0){6}}\put(46.5,22.5){\vector(1,0){19}}
\end{picture}}
\end{center}
\caption{Construction of an associated permutation}
\end{figure}

For a subset $A\subseteq[n]$, the complement set $\overline{A}$ is
obtained by replacing each $i\in A$ by $\overline{i}:=n+1-i$. The
\emph{reverse path}  of $w$ is the path $\overline{w}$ whose
$i$-th step is North (resp. East, Null, South-East) if and only if
the $(n+1-i)$-th step of $w$ is South-East (resp. East, Null,
North). Clearly, if $w\in \Omega_n^k$ with
$\lambda(w)=(\O,\C,\S,\T)$, then
 we can also define $\overline{w}$
as the unique path satisfying
$\type(\overline{w})=\overline{\type(w)}:=(\overline{\C},\overline{\O},\overline{\S},
\overline{\T})$.

\begin{lemma}\label{lemma:reverse}
The mapping $w\mapsto \bar w$ is an involution on $\Omega_n^k$.
Moreover,
\begin{enumerate}
\item  For $i\in [n]$, we have
$\h_i(\overline{w})=y_{\overline{i}+1}(w)$. In particular, we
have:
\begin{align}\label{eq:height}
\h_i(\overline{w})=\left\{%
\begin{array}{ll}
   \h_{\overline{i}}(w)-1, & \hbox{if $i\in\O(\overline{w})$;} \\
\h_{\overline{i}}(w), & \hbox{if
$i\in\T\cup\S(\overline{w})$;} \\
\h_{\overline{i}}(w)+1, & \hbox{if
$i\in\C(\overline{w})$.} \\
\end{array}%
\right.
\end{align}
\item Suppose $|\O(w)|=r$ and let $\sigma$ and $\sigma'$ be the associated permutations of
$w$ and $\overline{w}$ respectively. Then for any~$j\in[r]$,
\begin{equation}\label{eq:permasso}
\sigma'(j)=r+1-{\sigma}^{-1}(r+1-j).
\end{equation}

\end{enumerate}
\end{lemma}
\begin{proof}

(1) By definition of $\overline{w}$, the height of the $i$-step of
$\bar w$ corresponds to the  height of the $\bar i+1$-th step of
$w$, so $\h_i(\overline{w})=\h_{\overline{i}+1}(w)$. Eq.
\eqref{eq:height} follows then from the fact that $\lambda(\bar
w)=\overline{\lambda(w)}$ and the equation:
$$\h_{i+1}(w)= \left\{
                      \begin{array}{ll}
                        \h_{i}(w)+1 & \hbox{if $i\in\O(w)$;} \\
                        \h_{i}(w) & \hbox{if $i\in\T\cup\S(w)$;} \\
                        \h_{i}(w)-1 & \hbox{if $i\in\C(w)$.}
                      \end{array}
                    \right.$$
(2)  By definition, for any $j\in [r]$, the mapping $\sigma$ maps
the $j$-th North step of $w$ to the $\sigma(j)$-th South-East step
of $w$. Equivalently the mapping $\sigma'$ maps the
$r+1-\sigma(j)$-th North step of $\bar w$ to the $r+1-j$-th
South-East step of $\bar w$. In other words, we have
$\sigma'(r+1-\sigma(j))=r+1-j$.
 Substituting $r+1-\sigma(j)$ by $i$
yields the desired result.
\end{proof}

We have now all the ingredients to  define our involution
$\varphi$ on $\Delta_n^k$.

\noindent{\bf Involution $\varphi$.} Let $h=(w,\gamma)\in
\Delta_n^k$ such that  $\O\cup\S(w)=\{i_1,i_2,\cdots,i_k\}_{<}$,
$\T(w)=\{t_1,t_2,\cdots,t_u\}_{<}$ and
$\C(w)=\{c_1,c_2,\cdots,c_r\}_{<}$. If
$\O\cup\S(\overline{w})=\{i'_1,i'_2,\cdots,i'_k\}_{<}$,
$\T(\overline{w})=\{t'_1,t'_2,\cdots,t'_u\}_{<}$ and
$\C(\overline{w})=\{c'_1,c'_2,\cdots,c'_r\}_{<}$, then
$\varphi(h)=(\overline{w},\xi)$, where $\xi=(\xi_i)_{1\leq i\leq
n}$ is defined as follows:
\begin{align}\label{eq:defvarphi}
\xi_i=\left\{%
\begin{array}{ll}
    \gamma_{i_m}, & \hbox{if $i=i_m'$ for $m\in[k]$;} \\
    \gamma_{t_{u+1-m}}, & \hbox{if $i=t_m'$ for $m\in[u]$;} \\
    \gamma_{c_{\sigma(r+1-m)}}, & \hbox{if $i=c_m'$ for $m\in[r]$;} \\
\end{array}%
\right.
\end{align}
where $\sigma$ is  the associated permutation of $w$.
\begin{exam}\label{exam:invopath}
Consider the  path diagram $h=(w,\gamma)$ in the Figure~4.
 It is easy to see that the permutation associated to $w$
is $\sigma=321$. It follows that
$\xi_{c'_1}=\gamma_{c_{\sigma(3)}}=\gamma_{c_{1}}$,
$\xi_{c'_2}=\gamma_{c_{\sigma(2)}}=\gamma_{c_{2}}$ and
$\xi_{c'_3}=\gamma_{c_{\sigma(1)}}=\gamma_{c_3}$. The image
$\varphi(h)$ of $h$ is given below.
\end{exam}
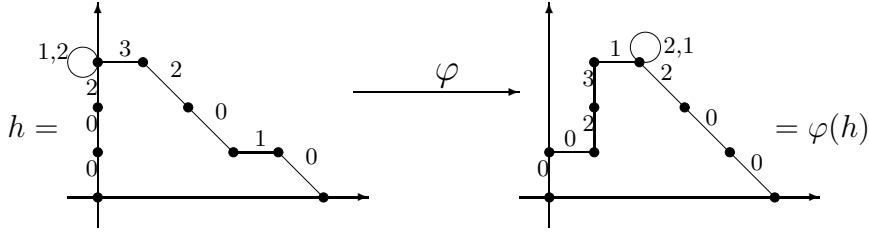
\begin{figure}[h]\
\begin{center}
{\setlength{\unitlength}{0.4mm}
\begin{picture}(230,70)(0,0)
\put(0,0){\circle*{3}}\put(0,15){\circle*{3}}\put(0,30){\circle*{3}}
\put(0,45){\circle*{3}}\put(15,45){\circle*{3}}
\put(30,30){\circle*{3}}\put(45,15){\circle*{3}}\put(60,15){\circle*{3}}
\put(75,0){\circle*{3}}
\put(-10,0){\vector(1,0){100}}\put(0,-10){\vector(0,1){75}}
\put(0,0){\linethickness{0.2mm}\line(0,1){45}}\put(0,45){\linethickness{0.2mm}\line(1,0){15}}
\put(15,45){\linethickness{0.2mm}\line(1,-1){30}}\put(45,15){\linethickness{0.2mm}\line(1,0){15}}
\put(60,15){\linethickness{0.2mm}\line(1,-1){15}}
\put(-5,45){\linethickness{0.2mm}\circle{10}}
\put(-4,7){\scriptsize0}\put(-4,22){\scriptsize0}\put(-4,34){\scriptsize2}
\put(-20,46){\scriptsize1,2}\put(7,47){\scriptsize3}\put(24,40){\scriptsize2}
\put(39,26){\scriptsize0}\put(52,17){\scriptsize1}\put(69,11){\scriptsize0}
\put(150,0){\circle*{3}}\put(150,15){\circle*{3}}\put(165,15){\circle*{3}}
\put(165,30){\circle*{3}}\put(165,45){\circle*{3}}
\put(180,45){\circle*{3}}\put(195,30){\circle*{3}}
\put(210,15){\circle*{3}}\put(225,0){\circle*{3}}
\put(150,0){\line(0,1){15}}\put(150,15){\line(1,0){15}}
\put(165,15){\line(0,1){15}}\put(165,30){\line(0,1){15}}
\put(165,45){\line(1,0){15}}\put(180,45){\line(1,-1){15}}\put(195,30){\line(1,-1){15}}
\put(210,15){\line(1,-1){15}}
\put(182,50){\linethickness{0.2mm}\circle{10}}\put(140,0){\vector(1,0){100}}
\put(150,-10){\vector(0,1){75}}
\put(146,7){\scriptsize0}\put(155,17){\scriptsize0}\put(161,22){\scriptsize2}
\put(161,37){\scriptsize3}\put(170,47){\scriptsize1}\put(188,48){\scriptsize2,1}
\put(187,39){\scriptsize2}\put(202,24){\scriptsize0}\put(217,9){\scriptsize0}
\put(-30,20){$h=$}\put(225,20){$=\varphi(h)$}
\put(85,35){\vector(1,0){55}}\put(112,39){\large \bf{$\varphi$}}

\end{picture}}
\end{center}
\caption{An example of the involution $\varphi$}
\end{figure}
\bigskip
We now present the main result of this section.
\begin{prop}\label{thm:F}
The mapping $\varphi: (w,\gamma)\mapsto (\bar w,\xi)$   is an
involution on $\Delta_n^k$ such that
\begin{equation}\label{eq:thmvarphi}
\sum_{i\in\O\cup\S(\overline{w})}\xi_i
=\sum_{i\in\O\cup\S(w)}\gamma_i \quad\text{and}\quad
\sum_{i\in\T\cup\C(\overline{w})}\xi_i=\sum_{i\in\T\cup\C(w)}\gamma_i.
\end{equation}
\end{prop}

\pf We first show that the mapping $\varphi$ is well defined. It
suffices to show that:
\begin{itemize}
\item[(a)] $0\leq \xi_{i'_m}\leq x_{i'_m}(\overline{w})+\h_{i'_m}(\overline{w})$
for $m\in[k]$,
\item[(b)] $0\leq \xi_{t'_m}\leq \h_{t'_m}(\overline{w})-1$ for
$m\in[u]$, and
\item[(c)] $0\leq \xi_{c'_m}\leq \h_{c'_m}(\overline{w})-1$ for
$m\in[r]$.\\
\end{itemize}
Since  $t'_m=\overline{t_{u+1-m}}$ and
$c'_m=\overline{o_{r+1-m}}$, we have
$$\overline{t'_m}=\overline{\overline{t_{u+1-m}}}=t_{u+1-m}\quad\textrm{
and}\quad
\overline{c'_m}=\overline{\overline{o_{r+1-m}}}=o_{r+1-m}.
$$
\begin{itemize}
\item[(i)]
For $m\in[k]$, Lemma~\ref{usefullemma}(ii) implies that
$x_{i'_m}(\overline{w})+\h_{i'_m}(\overline{w})
=x_{i_m}(w)+\h_{i_m}(w)$ because both sides are equal to $m-1$. We have
(a) by invoking  $\xi_{i'_{m}}=\gamma_{i_{m}}$ and
 $0\leq \gamma_{i_{m}}\leq
x_{i_m}(w)+\h_{i_m}(w)$.

\item[(ii)] For $m\in[u]$, Lemma~\ref{lemma:reverse}(2) implies that
$\h_{t'_m}(\overline{w})=\h_{\overline{t'_m}}(w)=\h_{t_{u+1-m}}(w)$.
We obtain
(b)  by invoking $\xi_{t'_m}=\gamma_{t_{u+1-m}}$,
  $0\leq
\gamma_{t_{u+1-m}}\leq \h_{t_{u+1-m}}(w)-1$
and $\h_{t_{u+1-m}}(w)=\h_{t'_m}(\overline{w})$.

\item[(iii)] For $m\in[r]$, Lemma~\ref{lemma:reverse}(2) implies
$\h_{c'_m}(\overline{w})=\h_{c_{\sigma(r+1-m)}}(w)$ because
$\h_{c'_m}(\overline{w})=\h_{\overline{c'_m}}(w)+1=\h_{o_{r+1-m}}(w)+1$,
 which is equal to $\h_{c_{\sigma(r+1-m)}}(w)$ by definition of $\sigma$.
We derive (c)  by invoking that $\xi_{c'_m}=\gamma_{c_{\sigma(r+1-m)}}$ and
$0\leq\gamma_{c_{\sigma(r+1-m)}}\leq
\h_{c_{\sigma(r+1-m)}}(w)-1$.

\end{itemize}

Let  $\sigma'$ be the associated permutation of
$\overline{w}$. Consider the chain of bijections:
$$
\varphi^2:\; (w,\gamma)\stackrel{\varphi}{\longmapsto}(\overline{w},\gamma)
\stackrel{\varphi}{\longmapsto}(\overline{\overline{w}},\mu)=(w,\mu)
$$
 where
$\mu=(\mu_i)_{1\leq i\leq n}$ is defined by
\begin{align*}
\mu_i=\left\{%
\begin{array}{ll}
    \xi_{i'_m}, & \hbox{if $i=i_m$ for $m\in[k]$;} \\
    \xi_{t'_{u+1-m}}, & \hbox{if $i=t_m$ for $m\in[u]$;} \\
    \xi_{c'_{\sigma'(r+1-m)}}, & \hbox{if $i=c_m$ for $m\in[r]$.} \\
\end{array}%
\right.
\end{align*}

Therefore
\begin{align*}
\mu_{i_m}=\xi_{i'_m}=\gamma_{i_m}\;\text{for
$m\in[k]$}\quad\text{and}\quad
\mu_{t_m}=\xi_{t'_{u+1-m}}=\gamma_{t_m}\;\text{for $m\in[u]$}.
\end{align*}
By \eqref{eq:permasso} we have
$\sigma(r+1-\sigma'(j))=r+1-j$ and
\begin{align*}
\mu_{c_m}=\xi_{c'_{\sigma'(r+1-m)}}=\gamma_{c_{\sigma(r+1-\sigma'(r+1-m))}}
=\gamma_{c_m}\quad\text{for
$m\in[r]$}.
\end{align*}
Hence $\mu=\gamma$ and $\varphi$ is an involution.
Finally \eqref{eq:thmvarphi} follows from
\eqref{eq:defvarphi}.
\qed


\section{Bijection $\Gamma_\sigma$ from $\P_n^k$ to $\P_n^k(\sigma)$}
The {\it Lehmer code} of a permutation $\sigma$ in $S_n$ is the
sequence $c(\sigma)=(c_1, \ldots, c_n)$ of non negative integers
where the integer $c_i$ is defined by
$$
c_i = \# \{ j > i, \sigma(j) < \sigma(i) \}.
$$
We can recover  the permutation $\sigma$ from its  code
$c(\sigma)$ because $\sigma(i)$ equals the $(c_i+1)$-th element in
$[n] - \{ \sigma(1),...,\sigma(i-1) \}$. Therefore the mapping $c$
which associates to each permutation of $S_n$ its Lehmer code is a
bijection from $S_n$ to  $[0,n-1]\times [0,n-2]\times ... \times
[0,1]\times [0]$. Moreover $\inv\sigma=c_1+\ldots +c_n$. For our
purpose we need to define the $d$-code of a permutation $\sigma\in
S_n$ by
$$
d(\sigma):=(d_1,\ldots, d_n)=(c_{\sigma^{-1}(1)},\ldots, c_{\sigma^{-1}(n)} ).
$$
That is, the coordinate $d_i$ is the number of entries $\sigma(j)$ smaller than and
to the right of
$i$ in the sequence $\sigma(1)\ldots \sigma(n)$. In other words, we have
$$
d_i=\#\{j>\sigma^{-1}(i); \sigma(j)<i\}.
$$
\begin{lemma}
The mapping $d$ which associates to each permutation of $S_n$ its d-code is a bijection
from $S_n$ to  $P_n = [0]\times [0,1]\times ... \times [0,n-2]\times [0, n-1]$. Moreover
$\inv\sigma=d_1+\ldots +d_n$.
\end{lemma}
For example, if
$\sigma = 8 \,6 \,3 \,4\, 7 \,5\, 2 \,1$, then
\begin{align*}
c(\sigma) = 7\, 5\, 2\, 2 \,3\, 2\, 1\, 0,\\
d(\sigma) = 0\, 1\, 2\, 2 \,2\, 5\, 3\, 7.
\end{align*}

\begin{lemma} \label{usefullemma2}
For any $(w,\gamma)\in \Delta_n^k$ with $\O\cup\S(w)=\{i_1,
i_2,\ldots,i_k\}_<$,
 the sequence $(\gamma_{i_{1}},\ldots, \gamma_{i_{k}})$ is the
 $d$-code of some permutation
$\sigma\in S_k$.
 \end{lemma}
 \begin{proof}
It suffices to verify that $0\leq \gamma_{i_{j}}\leq j-1$ for
$j=1,\ldots, k$, but this is obvious in view of
Lemma~\ref{usefullemma}(ii).
\end{proof}

For any permutation $\sigma\in S_k$, let $\Delta_n^k(\sigma)$ be
the set of path diagrams $(w,\gamma)\in \Delta_n^k$ such that
$d(\sigma)=(\gamma_{i_{1}},\ldots, \gamma_{i_{k}})$, where
$\O\cup\S(w)=\{i_1, i_2,\ldots,i_k\}_<$.

\begin{lemma} \label{sigmapath-permutation}
For any $\sigma\in S_k$ we have
\begin{itemize}
\item[(i)] the restriction of $\varphi$ on $\Delta_n^k(\sigma)$ is
an involution;
\item[(ii)] the restriction of $\Phi$ on $\Delta_n^k(\sigma)$ is a
bijection from $\Delta_n^k(\sigma)$ to $\P_n^k(\sigma)$.
\end{itemize}
\end{lemma}
\begin{proof}
Let $(w,\gamma)\in \Delta_n^k(\sigma)$ with
$\O\cup\S(w)=\{i_1,i_2,\ldots, i_k\}_\leq$. Then
$(\gamma_{i_{1}},\ldots, \gamma_{i_{k}})$ is the $d$-code of
$\sigma$. (i) follows directly from the definition of $\varphi$.
Let $\Phi(w,\gamma)=\pi$. Suppose
$\pi=B_{\tau(1)}/B_{\tau(2)}/\ldots /B_{\tau(k)}$
 is a $\tau$-partition in $\OP_n^k$ and $d(\tau)=(d_1,\ldots, d_k)$. Hence
 $B_1/B_2/\ldots /B_k$ is a partition in $\P_n^k$ and $\min B_j=i_j$ for $j\in [k]$.
By definition of $\Phi$, in the $i_{j}$-th step, we create a new block $B_{j}$
 with  $i_{j}$ in $T_{i_j}$ so that
 there are $\gamma_{i_{j}}$ blocks to the right of
 the  block $B_{j}$, i.e., $d_{j}=\gamma_{i_{j}}$. Namely $\tau=\sigma$.

Conversely, given $\pi\in \P_n^k(\sigma)$, then there is a path
diagram $(w,\gamma)\in \Delta_n^k(\tau)$, for some $\tau\in S_k$,
such that $\Phi(w,\gamma)=\pi$. As $\Phi$ is a bijection, we get
 $d(\sigma)=d(\tau)$, so $\sigma=\tau$. This completes the proof
 of  (ii).
 \end{proof}

It is convenient to introduce the following abbreviations:
\begin{align}\label{abbreviation}
\cls:=\lcs+\rcs, \quad \opb:=\lob+\rob\quad\textrm{and}\quad \sbg:=\lsb+\rsb.
\end{align}

\begin{prop}\label{prop:propriétés du type}
For any $\pi\in\OP_n^k$,
\begin{align}
\cls(\pi)&= \sum_{i\in\, \clos(\pi)}(n-i),\label{eq1}\\
\opb(\pi)&= \sum_{i\in\, \open(\pi)}(i-1),\label{eq2}\\
\sbg(\pi)&= \sum_{i\in\, \C(\pi)}i-\sum_{i\in\,
O(\pi)}i+k-n.\label{eq3'}
\end{align}
\end{prop}

\pf First of all, equations \eqref{eq1} and \eqref{eq2} follow
immediately from the fact that for any $i\in [n]$ the number of
integers greater (resp. smaller) than $i$ in $[n]$ equals $n-i$
(resp. $i-1$). Indeed, by definition \eqref{abbreviation}, the
statistic $\cls(\pi)$ amounts to count, for each closer $i$ of
$\pi$, the number of  the integers greater than $i$ in $[n]$.
Similarly we get \eqref{eq2}. Now, suppose $\pi=B_1/\cdots/B_k$,
then
\begin{align*}
\sb(\pi)&=\sum_{i=1}^k|\{j\,:\min(B_i)<j<\max(B_i),\,j\notin B_i\}|\\
       &=\sum_{i=1}^k(\max(B_i)-\min(B_i)+1-|B_i|)\\
       &=\sum_{i=1}^k\max(B_i)-\sum_{i=1}^k\min(B_i)+k-n,
\end{align*}
which is exactly \eqref{eq3'}.
\qed
\medskip

We are now ready to construct   a bijection from (no ordered) set
partitions to ordered set partitions and state the main theorem of
this section.
\begin{thm}\label{thm:sigmainvolution} For any $\sigma\in S_k$
there is a bijection $\Gamma_\sigma : \P_n^k\to \P_n^k(\sigma)$
such that for any $\pi\in\P_n^k$,
\begin{itemize}
\item[(i)] $\lambda(\Gamma_\sigma(\pi))=\lambda(\pi)$;
\item[(ii)] $(\cls, \opb,\sb)\Gamma_\sigma(\pi)=(\cls,\opb,\sb)\pi$;
\item[(iii)] $\rsb_\tc \Gamma_\sigma(\pi)=\rsb_\tc\pi$.
\end{itemize}
\end{thm}
\begin{proof}
Let $\sigma, \tau$ be  two permutations in $S_k$ with
$d(\sigma)=(d_1,\ldots, d_k)$ and $d(\tau)=(d_1',\ldots, d_k')$.
For any path diagram $(w, \gamma)\in \Delta_n^k(\tau)$ we can
define a path diagram $g_{\tau,\sigma}(w, \gamma)=(w, \gamma')\in
\Delta_n^k(\sigma)$ as follows:
\begin{align}\label{eq:defgsigma}
\gamma_i'=\left\{%
\begin{array}{ll}
  d_j, & \hbox{if $i=i_j\in \O\cup \S(w)$;} \\
  \gamma_i, & \hbox{if $i\in\T\cup\C(w)$.}\\
\end{array}
\right.
\end{align}
Clearly the mapping $g_{\sigma,\tau}: \Delta_n^k(\sigma)\to
\Delta_n^k(\tau)$ is a bijection because
$g_{\sigma,\tau}^{-1}=g_{\tau,\sigma}$.

In particular, taking $\tau=\epsilon$ (the identity permutation),
then $g_\sigma:=g_{\epsilon,\sigma}$ is a bijection from
$\Delta_n^k(\epsilon)$ to $\Delta_n^k(\sigma)$. It follows that
\begin{align}\label{eq:defgamma}
\Gamma_\sigma:=\Phi\circ g_\sigma\circ \Phi^{-1}
\end{align}
is a bijection from $\P_n^k$ to $\P_n^k(\sigma)$. For any  $\pi\in
\P_n^k$, let $\Gamma_\sigma(\pi)=\pi'$.
The composition of mappings
is better understood by the following diagram:
$$
\Gamma_\sigma:\;\pi\stackrel{\Phi^{-1}}{\longmapsto}(w,
\gamma)\stackrel{g_\sigma}{\longmapsto}(w, \gamma')
\stackrel{\Phi}{\longmapsto}\pi'.
$$
By  definition of $\Phi$ and $g_\sigma$, we have
$\type(\pi')=\type(w)=\type(\pi)$, namely (i). (ii) follows from
Proposition~\ref{prop:propriétés du type}. Combining
\eqref{eq:thmPsi} and \eqref{eq:defgsigma} yields (iii):
$$
{\rsb}_{\tc}\,\pi'=\sum_{i\in\T\cup\C(w)}\gamma_i'=\sum_{i\in\T\cup\C(w)}\gamma_i={\rsb}_{\tc}\,\pi.
$$

This completes the proof of Theorem~~\ref{thm:sigmainvolution}.
\end{proof}

\begin{exam} Let $\pi=1\;5\;7/2\;4\;10/3\;8/6/9\in \P_{10}^5$
and $\sigma=4\,3\,1\,5\,2$. Hence $opener(\pi)=\{1,2,3,6,9\}_{<}$
 and $d(\sigma)=0\,0\,2\,3\,1$. The corresponding path diagrams
$(w, \gamma)$  and $(w, \gamma')$ are given by
\begin{align*}
w&=(0,0),\, (0,1),\,(0,2),\,(0,3),\,(0,3),\,(0,3),\,(1,3),\,(2,2),\,(3,1),\,(4,1),\,(5,0),\\
\gamma&=(\mathbf{0},\,\mathbf{0},\,\mathbf{0},\, 1,\,2,\,\mathbf{0},\,2,\, 0,\, \mathbf{0},\, 0 ),\\
\gamma'&=(\mathbf{0},\,\mathbf{0},\,\mathbf{2},\, 1,\,2,\,\mathbf{3},\,2,\, 0,\, \mathbf{1},\, 0 ).
\end{align*}

 The mapping $g_\sigma: h\mapsto h'$ is illustrated in Figure~5.
\begin{figure}[h]
$$
{\setlength{\unitlength}{0.4mm}
\begin{picture}(180,70)(0,0)
\put(-80,30){\makebox(3,4){$h\;=$}}
\put(-50,0){\circle*{3}}\put(-50,15){\circle*{3}}
\put(-50,30){\circle*{3}}\put(-50,45){\circle*{3}}\put(-35,45){\circle*{3}}
\put(-20,30){\circle*{3}}\put(-5,15){\circle*{3}}\put(10,15){\circle*{3}}\put(25,0){\circle*{3}}
\put(-60,0){\vector(1,0){100}}\put(-50,-5){\vector(0,1){65}}
\put(-50,0){\linethickness{0.2mm}\line(0,1){45}}
\put(-50,45){\linethickness{0.2mm}\line(1,0){15}}
\put(-35,45){\linethickness{0.2mm}\line(1,-1){30}}
\put(-5,15){\linethickness{0.2mm}\line(1,0){15}}
\put(10,15){\linethickness{0.2mm}\line(1,-1){15}}
\put(-55,45){\linethickness{0.2mm}\circle{10}}

\put(-54,7){\scriptsize$\bf{0}$}\put(-54,22){\scriptsize$\bf{0}$}\put(-54,34){\scriptsize$\bf{0}$}
\put(-72,46){\scriptsize$1,2$}\put(-43,47){\scriptsize$\bf{0}$}\put(-26,40){\scriptsize$2$}
\put(-11,26){\scriptsize$0$}\put(2,17){\scriptsize$\bf{0}$}\put(19,9){\scriptsize$0$}
\put(36,27){$\xrightarrow[\phantom{2cm}]{g_\sigma}$}

\put(70,30){\makebox(3,4){$h'\;=$}}
\put(100,0){\circle*{3}}\put(100,15){\circle*{3}}\put(100,30){\circle*{3}}
\put(100,45){\circle*{3}}\put(115,45){\circle*{3}}
\put(130,30){\circle*{3}}\put(145,15){\circle*{3}}\put(160,15){\circle*{3}}\put(175,0){\circle*{3}}
\put(90,0){\vector(1,0){100}}\put(100,-5){\vector(0,1){63}}
\put(100,0){\linethickness{0.2mm}\line(0,1){45}}\put(100,45){\linethickness{0.2mm}\line(1,0){15}}
\put(115,45){\linethickness{0.2mm}\line(1,-1){30}}\put(145,15){\linethickness{0.2mm}\line(1,0){15}}
\put(160,15){\linethickness{0.2mm}\line(1,-1){15}}
\put(95,45){\linethickness{0.2mm}\circle{10}}

\put(96,7){\scriptsize$\bf{0}$}\put(96,22){\scriptsize$\bf{0}$}\put(96,34){\scriptsize$\bf{2}$}
\put(78,46){\scriptsize$1,2$}\put(107,47){\scriptsize$\bf{3}$}\put(124,40){\scriptsize$2$}
\put(139,26){\scriptsize$0$}\put(152,17){\scriptsize$\bf{1}$}\put(169,9){\scriptsize$0$}
\end{picture}}
$$
\caption{Mapping $g_\sigma$}
\end{figure}
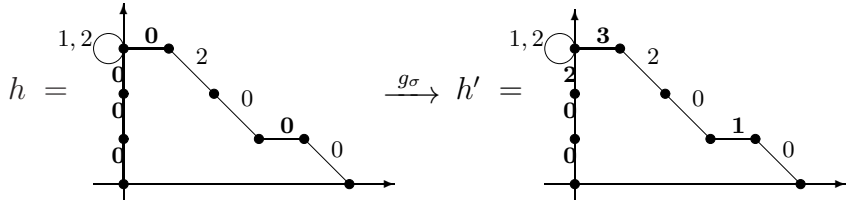

Finally, we get
$\Gamma_{\sigma}(\pi)=\Phi(h')=6/3\;5\;7/1\;4\;10/9/2\;8\in
\P_{10}^5(\sigma)$. Note that $\pi'$ is not a rearrangement of the
blocks in $\pi$.

\end{exam}

\section{Proof of Theorem~\ref{thm:sym}}
Consider the mapping
$$
\Xi=\Phi \circ \varphi\circ \Phi^{-1}:
\OP_n^k\longrightarrow \OP_n^k.
$$

For example, if $\pi=6/3\;5\;7/1\;4\;10/9/2\;8$, then
$\Phi^{-1}(\pi)=h$ is given in  Figure~2, while $\varphi(h)$ is
given in Figure~4. Finally we get
$\Xi(\pi)=\Phi(\varphi(h))=4\;6\;8/3\;7\;10/1\;9/5/2$.

Clearly
 the mapping
 $\Xi$ is an involution. For any  fixed  $\sigma\in\S_k$,
 Lemma~\ref{sigmapath-permutation} implies that
the restriction of $\Xi$ on $\P_n^k(\sigma)$ is stable.

 For any $\pi\in\P_n^k(\sigma)$, let $\Xi(\pi)=\pi'\in \P_n^k(\sigma)$.
Suppose that  $h=(w,\gamma)=\Phi^{-1}(\pi)$, $h'=(\bar
w,\xi)=\varphi(h)$ and $\pi'=\Phi(h')$. Then
Lemmas~\ref{usefullemma} and \ref{sigmapath-permutation} imply
that
\begin{align}\label{eq:inv=maj}
\Inv\pi=\inv\sigma=\Inv\pi'.
\end{align}
By Theorem~\ref{thmPhi} we see that
$\type(\pi')=\type(\overline{ w})$, $\type(\pi)=\type(w)$ and
$$
{\rsb}_{\tc}\,\pi'=\sum_{i\in\T\cup\C(\overline{w})} \xi_i
\quad\text{and}\quad {\rsb}_{\tc}\,\pi=\sum_{i\in\T\cup\C(w)}
\gamma_i.
$$
It follows from \eqref{eq:thmvarphi} that
\begin{align}\label{eq:rst}
{\rsb}_{\tc}\;\pi'={\rsb}_{\tc}\;\pi.
\end{align}
As $\type(\overline{ w})=\overline{\type(w)}$, by
Proposition~\ref{prop:propriétés du type},  we obtain
\begin{align}\label{eq:point}
(\cls,\opb,\sbg)\pi'=(\opb,\cls,\sbg)\pi.
\end{align}
Furthermore, on $\OP_n^k$, the following equations  hold true (see
\cite[Lemma~4.6]{IKZ}):
\begin{align}
\mak+\bInv&=\cls+\rsb_\tc+\Inv,\nonumber\\
\mak'+\bInv&=\opb+\rsb_\tc+\Inv,\label{eq:call}\\
\cinvLSB&=k(k-1)+\sb-\rsb_\tc-\Inv.\nonumber
\end{align}
It follows  from \eqref{eq:rst}, \eqref{eq:point} and
 \eqref{eq:call} that
\begin{align*}
&(\MAK+\bInv,\MAK'+\bInv,\cinvLSB)\,\pi'\\
&=(\MAK'+\bInv,\MAK+\bInv,\cinvLSB)\,\pi.
\end{align*}
This completes the first part of Theorem~\ref{thm:sym}.

In view of \eqref{eq:call} it is easy to see that \eqref{eq:KZWWbis1General} is equivalent to the following two identities:
\begin{equation}\label{eq3}
\sum_{\pi\in\,\P_n^k(\sigma)}\,
p^{\cls\,\pi+{\rsb}_{\tc}\,\pi}\,q^{\sbg\,\pi-{\rsb}_{\tc}\,\pi}=
\sum_{\pi\in\,\P_n^k(\sigma)}\,
p^{\opb\,\pi+{\rsb}_{\tc}\,\pi}\,q^{\sbg\,\pi-{\rsb}_{\tc}\,\pi},
\end{equation}
and
\begin{align}\label{eq4}
\sum_{\pi\in\,\P_n^k(\sigma)}\,
p^{\opb\,\pi+{\rsb}_{\tc}\,\pi}\,q^{\sbg\,\pi-{\rsb}_{\tc}\,\pi}
=S_{p,q}(n,k).
\end{align}
Now, \eqref{eq3} follows immediately from \eqref{eq:rst} and
 \eqref{eq:point}.
According to Theorem~\ref{thm:sigmainvolution}, it suffices to
prove \eqref{eq4} in the case of $\sigma=\varepsilon$.
But, on $\P_n^k=\P_n^k(\varepsilon)$,
since $\lob={\rsb}_{\os}=0$,
there hold
\begin{align*}
\opb+{\rsb}_{\tc}&=\rob+\rsb=\rcb,\\
\sbg-{\rsb}_{\tc}&=\lsb=\lsb+\lob=\lcb.
\end{align*}
Hence
the left-hand side of \eqref{eq4} is equal to
\begin{align}\label{eq5}
\sum_{\pi\in\,\P_n^k}\,
p^{\opb\,\pi+{\rsb}_{\tc}\,\pi}\,q^{\sbg\,\pi-{\rsb}_{\tc}\,\pi}
&=\sum_{\pi\in\,\P_n^k}\,
p^{\rcb\,\pi}\,q^{\lcb\,\pi},
\end{align}
and  \eqref{eq4} follows   by applying \eqref{eq:WW}.
The proof of Theorem~\ref{thm:sym} is thus completed.

\begin{rmk} For any mahonian permutation statistic $\mah$ we can
define a statistic $\Mah$ on ordered set partitions by
$\Mah\pi=\mah \sigma$ if $\pi\in \P_n^k(\sigma)$ and derive from
the above proof the following identities:
\begin{align}\label{eq:5.8}
\begin{gathered}
\sum_{\pi\in\,\OP_n^k}
p^{\cls\,\pi+{\rsb}_{\tc}\,\pi}\,q^{\sbg\,\pi-{\rsb}_{\tc}\,\pi}\,
t^{\Mah\,\pi} = \; [k]_t!\,S_{p,q}(n,k),\\
\sum_{\pi\in\,\OP_n^k}
p^{\opb\,\pi+{\rsb}_{\tc}\,\pi}\,q^{\sbg\,\pi-{\rsb}_{\tc}\,\pi}\,
t^{\Mah\,\pi} = \; [k]_t!\,S_{p,q}(n,k).
\end{gathered}
\end{align}
In particular, by taking $\mah=\inv$ we recover essentially
Theorem~\ref{thm:IKZ}, while by taking $\mah=\maj$ we derive that
the statistics
\begin{align*}
\cls+ {\rsb}_{\tc}+\Maj,\quad \opb+{\rsb}_{\tc}+\Maj, \quad
k(k-1)+\sb-\rsb_\tc-\Maj
\end{align*}
are Euler-Mahonian.
\end{rmk}
\begin{Remark} For ordinary partitions,
there is a similar bijection,  simpler than $\Xi$,  using Motzkin paths.
We sketch this bijection  below.

A {\it Motzkin path} of length $n$ is a lattice path in the plane of
integer lattice $\mathbb
Z^2$ from $(0,0)$ to $(n,0)$, consisting of NE-steps $(1,1)$, E-steps $(1,0)$,
and SE-steps $(1,-1)$, which never passes below the
$x$-axis. Let $\mathcal{D}_n$ be the set
of Motzin path diagrams $(\omega, \gamma)$,  where
$\omega$ is a Motzkin path of length $n$ and $\gamma=(\gamma_1,\ldots, \gamma_n)$ is a
sequence of labels such that if the $i$-th step is NE, then $\gamma_i=1$,
if the $i$-th step is SE, then $1\leq \gamma_i\leq h_i$,
and if the $i$-th step is E, then $1\leq \gamma_i\leq h_i+1$,
where $h_i$ is the height of the $i$-th step.

For each partition $\pi$ in $\P_n$ we can construct its
traces $\{T_i\}_{0\leq i\leq n}$.
Let $\gamma_i-1$ be the number of incomplete blocks to the left of
the block containing $i$ in $T_i$. As shown in \cite{KsZe}
we can construct a bijection $f: \pi\longmapsto (w, \gamma)$ from $\P_n$ to $\mathcal{D}_n$
 as follows:
 For $i=1,\ldots, n$, if
$i\in \O(\pi)$ (resp. $\S(\pi)$) we draw a NE-step (resp. E-step) with label
$1$ (resp. $\gamma_i$) and
if $i\in \C(\pi)$ (resp. $\T(\pi)$) we draw a SE-step
(resp. E-step) with label $\gamma_i$.

An example is given in Figure~6.

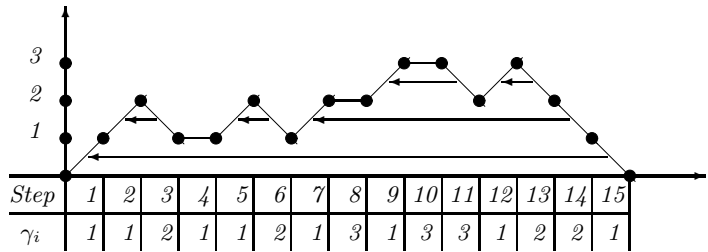
\begin{figure}[h]
\begin{center}
{\setlength{\unitlength}{0.5mm}
\begin{picture}(190,70)(0,-30)
\put(0,0){\circle*{3}}\put(-15,0){\linethickness{0.3mm}\vector(1,0){185}}
\put(-15,-10){\linethickness{0.3mm}\line(1,0){165}}
\put(0,-20){\linethickness{0.3mm}\vector(0,1){65}}
\put(0,0){\linethickness{0.2mm}\line(1,1){11.5}}
\put(10,10){\circle*{3}}
\put(10,10){\linethickness{2mm}\line(1,1){11.5}}
\put(20,20){\circle*{3}}
\put(20,20){\linethickness{2mm}\line(1,-1){11.5}}
\put(30,10){\circle*{3}}
\put(30,10){\linethickness{0.2mm}\line(1,0){10}}
\put(40,10){\circle*{3}}
\put(40,10){\linethickness{2mm}\line(1,1){11.5}}
\put(50,20){\circle*{3}}
\put(50,20){\linethickness{0.2mm}\line(1,-1){11.5}}
\put(60,10){\circle*{3}}
\put(60,10){\linethickness{0.2mm}\line(1,1){11.5}}
\put(70,20){\circle*{3}}
\put(70,20){\linethickness{0.2mm}\line(1,0){10}}
\put(80,20){\circle*{3}}
\put(80,20){\linethickness{0.2mm}\line(1,1){11.5}}
\put(90,30){\circle*{3}}
\put(90,30){\linethickness{0.2mm}\line(1,0){10}}
\put(100,30){\circle*{3}}
\put(100,30){\linethickness{0.2mm}\line(1,-1){11.5}}
\put(110,20){\circle*{3}}
\put(110,20){\linethickness{0.2mm}\line(1,1){11.5}}
\put(120,30){\circle*{3}}
\put(120,30){\linethickness{0.2mm}\line(1,-1){11.5}}
\put(130,20){\circle*{3}}
\put(130,20){\linethickness{0.2mm}\line(1,-1){11.5}}
\put(140,10){\circle*{3}}
\put(140,10){\linethickness{0.2mm}\line(1,-1){11.5}}
\put(150,0){\circle*{3}}
\put(10,0){\linethickness{0.2mm}\line(0,-1){20}}
\put(20,0){\linethickness{0.2mm}\line(0,-1){20}}
\put(30,0){\linethickness{0.2mm}\line(0,-1){20}}
\put(40,0){\linethickness{0.2mm}\line(0,-1){20}}
\put(50,0){\linethickness{0.2mm}\line(0,-1){20}}
\put(60,0){\linethickness{0.2mm}\line(0,-1){20}}
\put(70,0){\linethickness{0.2mm}\line(0,-1){20}}
\put(80,0){\linethickness{0.2mm}\line(0,-1){20}}
\put(90,0){\linethickness{0.2mm}\line(0,-1){20}}
\put(100,0){\linethickness{0.2mm}\line(0,-1){20}}
\put(110,0){\linethickness{0.2mm}\line(0,-1){20}}
\put(120,0){\linethickness{0.2mm}\line(0,-1){20}}
\put(130,0){\linethickness{0.2mm}\line(0,-1){20}}
\put(140,0){\linethickness{0.2mm}\line(0,-1){20}}
\put(150,0){\linethickness{0.2mm}\line(0,-1){20}}
\put(-10,30){\scriptsize3}
\put(-10,20){\scriptsize2}
\put(-10,10){\scriptsize1}
\put(0,30){\circle*{3}}
\put(0,20){\circle*{3}}
\put(0,10){\circle*{3}}
\put(-15,-7){\scriptsize \textrm{Step}}
\put(-12,-17){\scriptsize $\gamma_i$}
\put(5,-7){\scriptsize 1}
\put(15,-7){\scriptsize 2}
\put(25,-7){\scriptsize 3}
\put(35,-7){\scriptsize 4}
\put(45,-7){\scriptsize 5}
\put(55,-7){\scriptsize 6}
\put(65,-7){\scriptsize 7}
\put(75,-7){\scriptsize 8}
\put(85,-7){\scriptsize 9}
\put(92,-7){\scriptsize 10}
\put(102,-7){\scriptsize 11}
\put(112,-7){\scriptsize 12}
\put(122,-7){\scriptsize 13}
\put(132,-7){\scriptsize 14}
\put(142,-7){\scriptsize 15}
\put(5,-17){\scriptsize 1}
\put(15,-17){\scriptsize 1}
\put(25,-17){\scriptsize 2}
\put(35,-17){\scriptsize 1}
\put(45,-17){\scriptsize 1}
\put(55,-17){\scriptsize 2}
\put(65,-17){\scriptsize 1}
\put(75,-17){\scriptsize 3}
\put(85,-17){\scriptsize 1}
\put(94,-17){\scriptsize 3}
\put(104,-17){\scriptsize 3}
\put(114,-17){\scriptsize 1}
\put(124,-17){\scriptsize 2}
\put(134,-17){\scriptsize 2}
\put(144,-17){\scriptsize 1}
\put(-20,-35){$\pi=\{1,4, 15\}/\{2,3\}/\{5,6\}/\{7,10, 13\}
/\{8\}/\{9,11\}/\{12,14\}$.}
\put(144,5){\linethickness{0.2mm}\vector(-1,0){138}}
\put(24,15){\linethickness{0.2mm}\vector(-1,0){8}}
\put(54,15){\linethickness{0.2mm}\vector(-1,0){8}}
\put(134,15){\linethickness{0.2mm}\vector(-1,0){68}}
\put(104,25){\linethickness{0.2mm}\vector(-1,0){18}}
\put(124,25){\linethickness{0.2mm}\vector(-1,0){8}}
\end{picture}}
\end{center}
\caption{A labeled Motzkin path of length 15 and the corresponding partition. }
\end{figure}

Define an involution $g: (w,\gamma)\mapsto (w',\gamma')$ on $\mathcal{D}_n$
 as follows:
First reverse the path $w$ by reading it from right to left, i.e.,
if $w=((i,y_i))_{0\leq i\leq n}$, then
$w'=((i, y_{n-i}))_{0\leq i\leq n}$, then pair
the NE-steps with SE-steps in $w$ two by two in the following way:
 each NE-step  at height $h$  corresponds to the first SE-step to its right at height $h+1$
(thus we establish a bijection between the SE-steps of $w$ and those of $w'$),
 attribute the label of  each SE-step of $w$ to the corresponding
 SE-step of $w'$, finally
the labels of NE-steps of $w'$ are 1 and the E-steps of $w'$ keep the same
label as in $w$.

Now, it is easy to see (cf. \cite{KsZe}) that the mapping
$$\Lambda=f^{-1}\circ g\circ f: \pi\stackrel{f}{\longmapsto}(w,\gamma)
\stackrel{g}{\longmapsto}(w',\gamma')\stackrel{f^{-1}}{\longmapsto}\pi'$$
is an involution on $\P_n^k$ such that
$$
\mak\,(\pi)=\rcb\,\Lambda(\pi)\quad \textrm{and}\quad \lcb\,(\Lambda(\pi))=\lcb\,(\pi).
$$

The involution applied to the example in Figure~1 is given in  Figure~2, where
$$
\mak\pi=\rcb\Lambda(\pi)=37,\quad \lcb\pi=\lcb\Lambda(\pi)=16.
$$
\begin{figure}[h]
\begin{center}
{\setlength{\unitlength}{0.5mm}
\begin{picture}(190,70)(0,-30)
\put(0,0){\circle*{3}}\put(-15,0){\linethickness{0.3mm}\vector(1,0){185}}
\put(-15,-10){\linethickness{0.3mm}\line(1,0){165}}
\put(0,-20){\linethickness{0.3mm}\vector(0,1){65}}
\put(0,0){\linethickness{0.2mm}\line(1,1){11.5}}
\put(10,10){\circle*{3}}
\put(10,10){\linethickness{2mm}\line(1,1){11.5}}
\put(20,20){\circle*{3}}
\put(20,20){\linethickness{0.2mm}\line(1,1){11.5}}
\put(30,30){\circle*{3}}
\put(30,30){\linethickness{2mm}\line(1,-1){11.5}}
\put(40,20){\circle*{3}}
\put(40,20){\linethickness{2mm}\line(1,1){11.5}}
\put(50,30){\circle*{3}}
\put(50,30){\linethickness{0.2mm}\line(1,0){10}}
\put(60,30){\circle*{3}}
\put(60,30){\linethickness{0.2mm}\line(1,-1){11.5}}
\put(70,20){\circle*{3}}
\put(70,20){\linethickness{0.2mm}\line(1,0){10}}
\put(80,20){\circle*{3}}
\put(80,20){\linethickness{0.2mm}\line(1,-1){11.5}}
\put(90,10){\circle*{3}}
\put(90,10){\linethickness{0.2mm}\line(1,1){11.5}}
\put(100,20){\circle*{3}}
\put(100,20){\linethickness{0.2mm}\line(1,-1){11.5}}
\put(110,10){\circle*{3}}
\put(110,10){\linethickness{0.2mm}\line(1,0){10}}
\put(120,10){\circle*{3}}
\put(120,10){\linethickness{0.2mm}\line(1,1){11.5}}
\put(130,20){\circle*{3}}
\put(130,20){\linethickness{0.2mm}\line(1,-1){11.5}}
\put(140,10){\circle*{3}}
\put(140,10){\linethickness{0.2mm}\line(1,-1){11.5}}
\put(150,0){\circle*{3}}
\put(10,0){\linethickness{0.2mm}\line(0,-1){20}}
\put(20,0){\linethickness{0.2mm}\line(0,-1){20}}
\put(30,0){\linethickness{0.2mm}\line(0,-1){20}}
\put(40,0){\linethickness{0.2mm}\line(0,-1){20}}
\put(50,0){\linethickness{0.2mm}\line(0,-1){20}}
\put(60,0){\linethickness{0.2mm}\line(0,-1){20}}
\put(70,0){\linethickness{0.2mm}\line(0,-1){20}}
\put(80,0){\linethickness{0.2mm}\line(0,-1){20}}
\put(90,0){\linethickness{0.2mm}\line(0,-1){20}}
\put(100,0){\linethickness{0.2mm}\line(0,-1){20}}
\put(110,0){\linethickness{0.2mm}\line(0,-1){20}}
\put(120,0){\linethickness{0.2mm}\line(0,-1){20}}
\put(130,0){\linethickness{0.2mm}\line(0,-1){20}}
\put(140,0){\linethickness{0.2mm}\line(0,-1){20}}
\put(150,0){\linethickness{0.2mm}\line(0,-1){20}}
\put(-10,30){\scriptsize3}
\put(-10,20){\scriptsize2}
\put(-10,10){\scriptsize1}
\put(0,30){\circle*{3}}
\put(0,20){\circle*{3}}
\put(0,10){\circle*{3}}
\put(-15,-7){\scriptsize \textrm{Step}}
\put(-12,-17){\scriptsize $\gamma_i$}
\put(5,-7){\scriptsize 1}
\put(15,-7){\scriptsize 2}
\put(25,-7){\scriptsize 3}
\put(35,-7){\scriptsize 4}
\put(45,-7){\scriptsize 5}
\put(55,-7){\scriptsize 6}
\put(65,-7){\scriptsize 7}
\put(75,-7){\scriptsize 8}
\put(85,-7){\scriptsize 9}
\put(92,-7){\scriptsize 10}
\put(102,-7){\scriptsize 11}
\put(112,-7){\scriptsize 12}
\put(122,-7){\scriptsize 13}
\put(132,-7){\scriptsize 14}
\put(142,-7){\scriptsize 15}
\put(5,-17){\scriptsize 1}
\put(15,-17){\scriptsize 1}
\put(25,-17){\scriptsize 1}
\put(35,-17){\scriptsize 2}
\put(45,-17){\scriptsize 1}
\put(55,-17){\scriptsize 2}
\put(65,-17){\scriptsize 3}
\put(75,-17){\scriptsize 3}
\put(85,-17){\scriptsize 2}
\put(94,-17){\scriptsize 1}
\put(104,-17){\scriptsize 2}
\put(114,-17){\scriptsize 1}
\put(124,-17){\scriptsize 1}
\put(134,-17){\scriptsize 2}
\put(144,-17){\scriptsize 1}
\put(-20,-35){$\Lambda(\pi)=\{1,12, 15\}/\{2,4\}/\{3,6,9\}/\{5,7\}/\{8\}/\{10,11\}/\{13,14\}$.}
\end{picture}}

\end{center}
\caption{The labeled Motzkin path and the associated  partition
corresponding to those of Figure~6. }
\end{figure}

Note that the mapping  $\Lambda=f^{-1}\circ g\circ f$ is a corrected version of
that given
in \cite{KsZe}.

\end{Remark}

\section{The second   path diagrams encoding $\Psi$ of ordered set partitions}

Recall that a block in a trace is a subset $B$ of
$[n]\cup\{\infty\}$ such that $B\cap [n]\neq\emptyset$. By
convention, the closer of $B$ is the greatest element of $B$.
Hence  $\max(B)=\infty$ if $\infty\in B$. Thus, the statistics
$\bDes$, $\bMaj$ and $\rsb_i$ can be easily extended to traces.

Let $T=B_1/\cdots/B_r$ be the  $(i-1)$-trace of a partition
($i\geq 1$). We can insert $i$ before $B_{1}$, between two
adjacent blocks $B_{j}$ and $B_{{j+1}}$, for $1\leq j\leq r-1$, or
after $B_{r}$. Label  these insertion positions  from left to
right by $0,1,\ldots,r$. We say that the position $j$ is
\emph{active} if $B_{j+1}$ is active.

Let $A$ and $D$ be the set of active positions and block descents in
$T$, respectively. We then label the $r+1$ positions in $T$ as
follows:
\begin{itemize}
\item label the right-most position by  $a_0=r$ and  the elements of $A\cup D$ from right to left by $a_1,a_2,\ldots,
a_t$. So $a_t<a_{t-1}<\cdots <a_2<a_1$,
\item label the remaining positions from left to right by
$a_{t+1},\ldots, a_r$. So $a_{t+1}<\cdots<a_r$.
\end{itemize}

\begin{lemma}\label{lemma:fundamental}
 Let $\ell$ be an integer satisfying $0\leq \ell\leq r$ and define $T'$ to
be the $i$-trace obtained by inserting in $T$ the block
$\{i,\infty\}$ or $\{i\}$ into position $a_{\ell}$. Then,
\begin{align}\label{eq:fonda}
\rsb_{i} T'+\bMaj T'-\bMaj T=\ell.
\end{align}
\end{lemma}

\begin{proof}
Clearly $A$ and $D$ are disjoint. Suppose
$A=\{o_1,o_2,\ldots,o_m\}_>$ and $D=\{d_1, d_2,\ldots,d_p\}_>$. So
$t=m+p$. Let $B$ be one of the two blocks $\{i,\infty\}$ or
$\{i\}$.

Since $T$ is a $(i-1)$-trace, the openers of $T$ are all smaller
than $i$. This implies that $\rsb_i T'$ is the number of active
blocks in $T'$ to the right of $B$. We distinguish four cases.
\begin{enumerate}

\item \emph{$\ell=0$}. Since $a_0=r$, we have
$T'=B_1/B_2/\ldots/B_r/B$ and
 the equation \eqref{eq:fonda} is obvious.

\item \emph{$1\leq \ell \leq t$} and
$a_{\ell}\in D$. Let $a_{\ell}=d_j$ for some $j$, $1\leq j\leq p$. The block
descent set of $T'$ is then
$$
\{d_p<\cdots<d_{j+1}<d_{j}+1<d_{j-1}+1<\cdots<d_1+1\}.
$$
Let $q$ be the greatest number such that $o_q > a_{\ell}$. Clearly,
$\rsb_i(T')=q$. We thus have that
$$
\rsb_{i} T'+\bMaj T'-\bMaj T=q+j.
$$
It suffices now to remark that $\ell=j+q$.\\

\item \emph{$1\leq \ell \leq t$} and $a_{\ell}\in A$. Then, $a_{\ell}=o_j$ for some $j$, $1\leq j\leq m$. Clearly, $\rsb_i(T')=j$.
Let $q$ be the greatest number such that $d_q > a_{\ell}$. Then, the
block descent set of $T'$ is
$$
\{d_p<\cdots<d_{q+1}<d_{q}+1<d_{q-1}+1<\cdots<d_1+1\}.
$$
  It then follows that
$$
\rsb_{i} T'+\bMaj T'-\bMaj T=j+q.
$$
It suffices now to remark that $\ell=j+q$.\\

\item \emph{$t+1\leq \ell\leq r$}. Note that $a_{\ell}$ is a position before
a complete block. Let $q$ and $s$ be the greatest numbers such that
$d_q>a_{\ell}$ and $o_s>a_{\ell}$. Clearly, $\rsb_i(T')=s$ and the block
descent set of $T'$ is
$$
\{d_p<\cdots<d_{q+1}<a_{\ell}+1<d_{q}+1<d_{q-1}+1<\cdots<d_1+1\}.
$$
It follows that
$$
\rsb_{i} T'+\bMaj T'-\bMaj T=s+q+a_{\ell}+1.
$$
It suffices now to remark that $a_{\ell}=\ell-(s+q)-1$. Indeed, $a_{\ell}$ is
equal to the number of positions in $T$ to the left of the position,
i.e. $(p-q)+(m-s)+(\ell-t-1)$, which is also equal to $=\ell-(q+s)-1$
since $t=p+m$.
\end{enumerate}
\end{proof}

We give an example to illustrate the above result.
\begin{exam} Let
$T=6\;11\,\infty/3\;5\;7/1\;4\;10\,\infty/9/2\;8$. So there are 6
insertion  positions and the set of  active positions and descents
is $A=\{0,2\}$ and $D=\{4\}$. Therefore
$$a_0=5,\; a_1=4,\; a_2=2,\; a_3=0,\quad\textrm{and}\quad
a_4=1,\;a_5=3.
$$

Denote by $T'$ the trace obtained by inserting in $T$ the block
$\{12\}$ into position $a_i$. Then we have: \small{
$$
\begin{array}{c|cccc}
i & a_i&&T' & \scriptsize{\rsb_{12}(T')+\bMaj(T')-\bMaj(T)}\\
\hline\\
0 &5 &&6\;11\,\infty/3\;5\;7/1\;4\;10\,\infty/9/2\;8/\mathbf{12}    & 0+4-4=0  \\
\hline\\
1 & 4&&6\;11\,\infty/3\;5\;7/1\;4\;10\,\infty/9/\mathbf{12}/2\;8   & 0+5-4=1\\
\hline\\
2 & 2&&6\;11\,\infty/3\;5\;7/\mathbf{12}/1\;4\;10\,\infty/9/2\;8    & 1+5-4=2\\
\hline\\
3 &0 &&\mathbf{12}/6\;11\,\infty/3\;5\;7/1\;4\;10\,\infty/9/2\;8    &2+5-4=3\\
\hline\\
4 & 1&&6\;11\,\infty/\mathbf{12}/3\;5\;7/1\;4\;10\,\infty/9/2\;8    & 1+7-4=4\\
\hline\\
5 & 3&&6\;11\,\infty/3\;5\;7/1\;4\;10\,\infty/\mathbf{12}/9/2\;8 &
0+9-4=5
\end{array}
$$
}
\end{exam}

We now construct a new bijection $\Psi:\Delta_n^k\mapsto\OP_n^k$
based on the above lemma.

\bigskip

\noindent{\bf Algorithm $\Psi$.} Let $h=(w,\gamma)\in\Delta_n^k$
be a path diagram. Set $T_0=\emptyset$. Construct recursively
$i$-skeletons $T_i$ for $i=1,\ldots, n$ such that $T_{i}$
 has $\h_{i+1}(w)$ active blocks and $x_{i+1}(w)$ complete blocks by the following process.
Suppose $T_{i-1}=B_{1}/B_{2}/\cdots /B_{\ell}$ and $T_{i-1}$ has
$\h_i(w)$ active blocks and $x_{i}(w)$ complete blocks. Label
the positions before $B_{1}$, between $B_{j}$ and $B_{{j+1}}$, for
$1\leq j\leq \ell-1$, and after $B_{\ell}$ from left to right by
$\{0,1,\cdots,\ell\}$. Extend $T_{i-1}$ to $T_{i}$ as
follows:
 \begin{itemize}
\item \emph{The $i$-th step of $w$ is North (resp. East)}:
 Let $A$ be the set of the  active positions in $T_{i-1}$
 and $D$ the set of block descents in $T_{i-1}$. Then, set $a_0=\ell$,
$\{a_1>a_2>\cdots>a_t\}=A\cup\D$ and let $a_{t+1}<\cdots<a_{\ell}$ be
the remaining positions. Then insert the block $\{i,\infty\}$
(resp. $\{i\}$) into position $a_{\gamma_i}$.

\item \emph{The $i$-th step of $w$ is Null (resp. South-East)}:
label the $\h_i(w)$ active blocks of $T_{i-1}$ from right to left
by $\{0,1,\cdots,\h_i(w)-1\}$. Then insert $i$ (resp. replace
$\infty$ by $i$) in the active block labeled by  $\gamma_i$.
\end{itemize}

Since $x_{n+1}(w)=k$ and $\h_{n+1}(w)=0$, it then follows that
$T_n$ has $k$ complete blocks and 0 active blocks, i.e.
$T_n\in\OP_n^k$. Define  $\Psi(h)=T_n$.

\begin{exam}\label{exam:Theta}
Let $h$ be the path diagram  in Figure~2,
 then
$$\Psi(h)= 6/3\;5\;7/9/1\;4\;10/2\;8.
$$
The step by step construction of $\Psi(h)$ is given in
Figure~\ref{fig:Theta}.
\end{exam}

\begin{figure}\label{fig:Theta}
\begin{center}
{\setlength{\unitlength}{0.4mm}
\begin{picture}(90,330)(0,-280)

\put(-30,30){\makebox(3,4){$h\;=$}}
\put(0,0){\circle*{3}}\put(0,15){\circle*{3}}\put(0,30){\circle*{3}}\put(0,45){\circle*{3}}\put(15,45){\circle*{3}}
\put(30,30){\circle*{3}}\put(45,15){\circle*{3}}\put(60,15){\circle*{3}}\put(75,0){\circle*{3}}
\put(-10,0){\vector(1,0){100}}\put(0,-10){\vector(0,1){75}}\put(0,0){\linethickness{0.2mm}\line(0,1){45}}\put(0,45){\linethickness{0.2mm}\line(1,0){15}}
\put(15,45){\linethickness{0.2mm}\line(1,-1){30}}\put(45,15){\linethickness{0.2mm}\line(1,0){15}}\put(60,15){\linethickness{0.2mm}\line(1,-1){15}}
\put(-5,45){\linethickness{0.2mm}\circle{10}}

\put(-4,7){\scriptsize0}\put(-4,22){\scriptsize0}\put(-4,34){\scriptsize2}
\put(-20,46){\scriptsize1,2}\put(7,47){\scriptsize3}\put(24,40){\scriptsize2}
\put(39,26){\scriptsize0}\put(52,17){\scriptsize1}\put(69,11){\scriptsize0}

\put(35,-5){\vector(0,-1){35}}\put(25,-29){\textbf{$\Psi$}}

\put(0,-250){\makebox(100,160){$ \scriptsize{
\begin{array}{cccccl}
i & \quad\textrm{step}_i\quad& \gamma_i &&& \quad T_i\\
\hline\\
0 &&  &&& \emptyset \\
\hline\\
1 & North& 0 &&& \textbf{1}\,\infty \\
\hline\\
2 & North& 0 &&& 1\,\infty/\textbf{2}\,\infty \\
\hline\\
3 & North& 2 &&& \textbf{3}\,\infty/1\,\infty /2\,\infty \\
\hline\\
4 & Null & 1 &&& 3\,\infty/1\;\textbf{4}\,\infty/2\,\infty \\
\hline\\
5 &Null & 2 &&& 3\;\textbf{5}\,\infty/1\;4\,\infty/2\,\infty \\
\hline\\
6 &East & 3 &&& \textbf{6}/3\;5\,\infty/1\;4\,\infty/2\,\infty \\
\hline\\
7 & South\text{-}East& 2 &&& 6/3\;5\;\textbf{7}/1\;4\,\infty/2\,\infty \\
\hline\\
8 & South\text{-}East& 0 &&& 6/3\;5\;7/1\;4\,\infty/2\;\textbf{8} \\
\hline\\
9 & East & 1 &&& 6/3\;5\;7/\textbf{9}/1\;4\,\infty/2\;8\\
\hline\\
10 & South\text{-}East & 0
&&&6/3\;5\;7/9/1\;4\;\textbf{10}/2\;8.\\
\hline\\
\end{array}}
$}}
\end{picture}}
\end{center}
\caption{the step by step construction of $\Psi(h)$}
\end{figure}

 To show that $\Psi$ is bijective, we give its inverse.

\noindent{\bf Algorithm $\Psi^{-1}$.} Let $\pi\in\OP_n^k$ be an
ordered partition. Let $w$ be the path defined by
$\type(w)=\type(\pi)$. Suppose
$\O\cup\S(w)=\{i_1<i_2<\cdots<i_k\}$. Let
$$\gamma_{i_j}=\rsb_{i_j}
\pi+ \bMaj T_{i_j}(\pi)-\bMaj T_{i_{j-1}}(\pi)\quad \textrm{
for}\quad i=1\cdots k$$
 and $\gamma_i=\rsb_i(\pi)$ if
$i\in\T\cup\S(w)$.  It is not hard to verify  that $(w,
(\gamma_i)_{1\leq i\leq n})$ is a walk diagram in $\Omega_n^k$ and
its image under $\Psi$ is $\pi$.

By definition of $\Psi$, for $1\leq i\leq n$,
 $\bMaj T_i(\pi)-\bMaj T_{i-1}(\pi)=0$ for any
$i\in\T\cup\C(\pi)$. Since $\bMaj \emptyset=0$ and $T_n(\pi)=\pi$,
splitting the set $[n]$ into $\O\cup \S(\pi)$ and $\T\cup \C(\pi)$
we see that
\begin{align}
\bMaj\,\pi =\sum_{i\in\O\cup\S(\pi)}\left(\bMaj T_i(\pi)-\bMaj
T_{i-1}(\pi)\right).
\end{align}
We summarize the main properties of $\Psi$ in  the following
theorem.
\begin{thm}\label{thm:Theta}
The mapping $\Psi:\Delta_n^k\mapsto\OP_n^k$ is a bijection such that
if $h=(w,\gamma)\in\Delta_n^k$ and $\pi=\Psi(h)$, then
$\type(\pi)=\type(w)$ and
\begin{align}
{\gamma}_i=\left\{%
\begin{array}{ll}
 \rsb_i(\pi)+\bMaj(T_i(\pi))-\bMaj(T_{i-1}(\pi)), & \hbox{if $i\in \O\cup\S(\pi)$;} \\
\rsb_i(\pi), & \hbox{if $i\in \T\cup\C(\pi)$.} \\
\end{array}%
\right.
\end{align}
Therefore
\begin{equation}\label{eq:thmTheta}
\sum_{i\in\O\cup\S(w)}\gamma_i={\rsb}_{\os}\,\pi+\bMaj \pi,
\quad\text{and}\quad
\sum_{i\in\T\cup\C(w)}\gamma_i={\rsb}_{\tc}\,\pi.
\end{equation}
\end{thm}

\section{Proof of Theorem~\ref{thm:refinesteinconj}}
Consider the mapping
$$\Upsilon:=\Psi \circ \Phi^{-1}:\OP_n^k\to
\OP_n^k.
$$

For example, if $\pi=6/3\;5\;7/1\;4\;10/9/2\;8$, then it follows
from Figure~2 and Figure~3 that $\Upsilon(\pi)=
6/3\;5\;7/9/1\;4\;10/2\;8$. Note that  $\Upsilon(\pi)$ is
generally not a rearrangement of the blocks of $\pi$.

The main properties of $\Upsilon$ is summarized in the following
lemma.

\begin{lemma}\label{thm:Upsilon}
The map $\Upsilon:\OP_n^k\mapsto\OP_n^k$ is a bijection such that
for any $\pi\in\OP_n^k$, we have:
\begin{itemize}
\item [(i)] $\type(\Upsilon(\pi))=\type(\pi)$,
\item [(ii)] ${\rsb}_{\tc}\;\Upsilon(\pi)={\rsb}_{\tc}\;\pi$,
\item [(iii)] $\MAJ \Upsilon(\pi)=\INV \pi$.
\end{itemize}
\end{lemma}
\begin{proof}

By definition the mapping $\Upsilon$ is a bijection. For
$\pi\in\OP_n^k$, let $h=(w,\gamma)=\Phi^{-1}(\pi)$ and
$\pi'=\Psi(h)$. Hence
 $\pi'=\Upsilon(\pi)$. By the construction of $\Psi$ and $\Phi$,
 it is clear  that  $\type(\pi')=\type(w)=\type(\pi)$.
Combining \eqref{eq:thmPsi}, \eqref{eq:thmTheta} and
\eqref{eq:defnalternative de INV}, we have
 \begin{align*}
\MAJ \pi'&={\rsb}_{\os}\,\pi+\bMaj
\pi'=\sum_{i\in\O\cup\S(w)}\gamma_i={\ros}_{\os}\,\pi=\INV \pi, \\
{\rsb}_{\tc}\,\pi'&=\sum_{i\in\T\cup\C(w)}\gamma_i={\rsb}_{\tc}\,\pi,
 \end{align*}
 completing the proof.
\end{proof}

Applying Lemma~\ref{thm:Upsilon} and
Proposition~\ref{prop:propriétés du type} we obtain
\begin{equation}\label{eq:7.1}
\begin{gathered}
(\cls+{\rsb}_{\tc},\opb+{\rsb}_{\tc}, \sbg-{\rsb}_{\tc},
\MAJ)\,\Upsilon(\pi)\\
=(\cls+{\rsb}_{\tc},\opb+{\rsb}_{\tc}, \sbg-{\rsb}_{\tc},
\INV)\,\pi.
\end{gathered}
\end{equation}
According to  \eqref{def:cinv-cmaj} and \eqref{def:mak-mak'}, the
following functional identities hold on $\OP_n^k$:
\begin{align*}
\MAK+\bMaj&=\lcs+\rcs+\rsb+\bMaj\\
          &=\cls+{\rsb}_{\tc}+{\rsb}_{\os}+\bMaj\\
&=\cls+{\rsb}_{\tc}+\MAJ, \\
\MAK'+\bMaj&=\lob+\rob+\rsb+\bMaj\\
           &=\opb+{\rsb}_{\tc}+{\rsb}_{\os}+\bMaj\\
&=\opb+{\rsb}_{\tc}+\MAJ,\\
\cmajLSB&=k(k-1)+\sbg-\rsb-\bMaj\\
        &=k(k-1)+\sbg-{\rsb}_{\tc}-{\rsb}_{\os}-\bMaj\\
&= k(k-1)+\sbg-{\rsb}_{\tc}-\MAJ.
\end{align*}
We thus derive from \eqref{eq:7.1}  and \eqref{eq:call} that
$$
\begin{gathered}
(\mak+\bMaj, \mak'+\bMaj, \cmajLSB)\pi'\\
= (\mak+\bInv, \mak'+\bInv, \cinvLSB)\pi.
\end{gathered}
$$
This completes the proof of Theorem~3.3.

\begin{rmk} From \eqref{eq:5.8} and \eqref{eq:7.1} we derive
 immediately the
following equivalent forms of  \eqref{eq:Mah1}:
\begin{align}\label{eq:MAJ-analog of IKZ 1}
\begin{gathered}
\sum_{\pi\in\,\OP_n^k}
p^{\cls\,\pi+{\rsb}_{\tc}\,\pi}\,q^{\sbg\,\pi-{\rsb}_{\tc}\,\pi}\,
t^{\MAJ \pi}= \; [k]_t!\,S_{p,q}(n,k),\\
\sum_{\pi\in\,\OP_n^k}p^{\opb\,\pi+{\rsb}_{\tc}\,\pi}\,q^{\sbg\,\pi-{\rsb}_{\tc}\,\pi}\,
t^{\MAJ \pi} = \; [k]_t!\,S_{p,q}(n,k).
\end{gathered}
\end{align}
\end{rmk}

\begin{rmk}
Composing $\Upsilon$ and $\Xi$ we obtain
 the  mapping
 $$\Theta=\Psi \circ \varphi\circ \Psi^{-1}: \OP_n^k\longrightarrow \OP_n^k.
 $$
  For example,
 if $\pi=6/3\;5\;7/9/1\;4\;10/2\;8$, then
$\Psi^{-1}(\pi)=h$, where  $h$ is the path diagram in
Figure~\ref{fig:Theta}. Therefore  $\varphi(h)$ is that given in
Example~\ref{exam:invopath}. The reader can verify that
$\Theta(\pi)=\Psi(\varphi(h))=4\;6\;8/1\;7\;10/3\;9/5/2$.

Obviuously,  the mapping $\Theta$ is an involution on $\OP_n^k$
satisfying, for any $\pi\in\OP_n^k$,
\begin{itemize}
\item [(i)] $\type(\Theta(\pi))=\overline{\type(\pi)}$,
\item [(ii)] ${\rsb}_{\tc}\;\Theta(\pi)={\rsb}_{\tc}\;\pi$,
\item [(iii)] $\MAJ \Theta(\pi)=\MAJ \pi$.
\end{itemize}
Therefore, by Proposition~\ref{prop:propriétés du type},
\begin{equation*}
(\opb+{\rsb}_{\tc},
\sbg-{\rsb}_{\tc},\MAJ)\,\Theta(\pi)=(\cls+{\rsb}_{\tc},
\sbg-{\rsb}_{\tc},\MAJ)\,\pi.
\end{equation*}
In other words, we have
\begin{equation}
\begin{gathered}
(\MAK+\bMaj,\MAK'+\bMaj,\cmajLSB)\,\Theta(\pi)\\
=(\MAK'+\bMaj,\MAK+\bMaj,\cmajLSB)\,\pi.
\end{gathered}
\end{equation}

\end{rmk}

\section{Proof of Theorem~\ref{thm:extMac}}
For any $\pi=B_1/B_2/\cdots /B_k$ in $\P_n^k$, we  have
$$
\R(\pi)=\{B_{\sigma(1)}/B_{\sigma(2)}/\cdots
/B_{\sigma(k)}|\sigma\in S_k\}.
$$
By \eqref{eq:inv,maj-perm-part} and \eqref{eq:defnalternative de
INV}, we have $\INV (B_{\sigma(1)}/B_{\sigma(2)}/\cdots
/B_{\sigma(k)})=\inv\sigma$. It follows by using Lehmer code (see
Section~6) that for any $\pi\in\P_n^k$,
\begin{equation}
\sum_{\pi\in\R(\pi)}q^{\INV \pi}=\sum_{\sigma\in S_k}q^{\inv
\sigma}=[k]_q!.
\end{equation}
It remains to show that
\begin{equation}
\sum_{\pi\in\R(\pi)}q^{\MAJ \pi}=[k]_q!.
\end{equation}

Let $C_k=\{(c_1,\ldots, c_k) :\;0\leq c_i\leq i-1\}$. We will
construct a  bijection $\beta_\pi:C_k\mapsto\R(\pi)$
 such that for any $\mathbf{c}=(c_1,\ldots, c_k)\in
C_k$, we have $\MAJ\beta_\pi(\mathbf{c})=\sum_{i=1}^k c_i$.

We construct recursively $i$-skeletons $T_i$, $1\leq i\leq n$,
such that $T_i$ has $\O(\pi)^{\leq i}-\C(\pi)^{\leq i}$ active
blocks and $(\C\cup\S(\pi))^{\leq i}$ complete blocks by the
following process.
 Set $T_0=\emptyset$ and suppose $T_{i-1}=B_{1}/B_{2}/\cdots /B_{l}$. Then $T_{i}$
is obtained from $T_{i-1}$ as follows:
 \begin{itemize}
\item $i\in\O(\pi)$ (resp. $\S(\pi)$): label the positions
before $B_{1}$, between $B_{j}$ and $B_{{j+1}}$, for $1\leq j\leq
l-1$, and after $B_{l}$ from left to right by $\{0,1,\cdots,l\}$.

 Let $F$ be the set of the positions before the active blocks in $T_{i-1}$
 and $D=\bDes(T_{i-1})$. Then, set $a_0=l$,
$\{a_1>a_2>\cdots>a_t\}=F\cup\D$ and let $a_{t+1}<\cdots<a_l$ be
the remaining positions. We then insert the block $\{i,\infty\}$
(resp. $\{i\}$) into position $a_{\gamma_i}$.

\item $i\in\T(\pi)$ (resp. $\C(\pi)$): insert $i$ (resp. replace $\infty$
by $i$) in the active block whose opener is the opener of the
block of $\pi$ which contains $i$.
\end{itemize}

 It is not difficult to see, via Lemma~\ref{lemma:fundamental}, that the above
procedure is well defined. Since $\O(P)^{\leq n}=\C(P)^{\leq n}$
and $(\C(P)\cup\S(P))^{\leq n}=|\C(P)\cup\S(P)|=k$, $T_n$ is a
$n$-skeleton with $0$ active blocks and $k$ complete blocks, i.e.
$T_n\in\OP_n^k$. Now, by construction, $T_n\in\R(\pi)$. We then
set $\beta_P(\mathbf{c})=T_n$.

To show that $\beta_P$ is bijective, we describe its inverse. Let
$\pi\in\R(\pi)$ and suppose
$$
\O\cup\S(\pi)=\{i_1<i_2<\cdots<i_k\}.
$$
 For $1\leq j\leq k$, let
$c_j=\rsb_{i_j} \pi+ \bMaj T_{i_j}(\pi)-\bMaj T_{i_{j-1}}(\pi)$.
It is then readily seen that $\beta_P(c_1\cdots c_k)=\pi$.

{\bf Remark.} A Foata style bijection which establishes directly
the equidistribution of $\INV$ and $\MAJ$ willl be given by the
first author in \cite{Ka}.
\section{Concluding remarks}

Recall that the $q$-Eulerian numbers $A_q(n,k)$ ($n\geq k\geq 0$)
of Carlitz~\cite{Ca2} are defined by
\begin{align}
A_q(n,k)=q^k[n-k]_qA_q(n-1,k-1) +[k+1]_qA_q(n-1,k),
\end{align}
 and have  the following combinatorial interpretation:
$$
A_q(n, k)=\sum_{\sigma}q^{\maj\, \sigma}
$$
where the summation is over all permutations of $[n]$  with $k$
descents.

The original motivation of  Steingr\'{\i}msson \cite{Stein} was
to give a direct combinatorial proof of the following
identity~\cite{ZZ}:
\begin{equation}\label{eq:zezh}
[k]_q!\,S_q(n,k)=\sum_{m=1}^{k}q^{k(k-m)}\,{n-m\brack
n-k}_q\,A_q(n,m-1).
\end{equation}
Though we have proved all the conjectures inspired by
\eqref{eq:zezh}, a direct combinatorial proof of \eqref{eq:zezh}
is still \emph{missing}. As proved in \cite{ZZ} the identity
\eqref{eq:zezh} is equivalent to Garsia's $q$-analogue of
Frobenius formula relating $q$-Eulerian numbers and $q$-Stirling
numbers of the second kind~(see \cite{Ga,GaRe}):
\begin{equation}\label{eq:qfrobenius}
\sum_{k=1}^n\frac{[k]_q!S_q(n,k)x^k}{(x;q)_{k+1}} =
\sum_{k=1}^\infty [k]_q^nx^k=\frac{\sum_{\sigma\in S_n} x^{1+\des
\sigma}q^{\maj \sigma}}{(x;q)_{n+1}}.
\end{equation}
Note that a combinatorial proof of  \eqref{eq:qfrobenius} has been given by Garsia and Remmel \cite{GaRe}.

Briggs and Remmel~\cite{BR} proved the following $p,q$-analogue of
Frobenius formula~\eqref{eq:qfrobenius}:
\begin{equation}\label{eq:pqfrobenius}
\sum_{k=1}^n\frac{[k]_{p,q}!{\hat S}_{p,q}(n,k) p^{{n-k+1\choose
2}+k(n-k)}x^k}{(xp^n;q/p)_{k+1}} =\frac{\sum_{\sigma\in
S_n}x^{\des \sigma+1}q^{\maj \sigma}p^{\comaj \sigma}}
{(xp^n;q/p)_{n+1}},
\end{equation}
where $\comaj\sigma=n\,\des\sigma-\maj\sigma$ and ${\hat
S}_{p,q}(n,k)$ is a  variante  of the $p,q$-Stirling numbers of
the second kind defined by the following recursion:
\begin{align}
{\hat S}_{p,q}(n,k)=q^{k-1}{\hat
S}_{p,q}(n-1,k-1)+p^{-n}[k]_{p,q}{\hat S}_{p,q}(n-1,k).
\end{align}

We would like to point out that \eqref{eq:pqfrobenius} and
 \eqref{eq:qfrobenius} are equivalent. Obviously \eqref{eq:qfrobenius} corresponds to the $p=1$ case of
\eqref{eq:pqfrobenius}. Conversely, since
$[k]_{q/p}!=p^{-{k\choose 2}}[k]_{p,q}!$ and
\begin{align*}
{\hat S}_{p,q}(n,k)&=p^{-{n-k+1\choose 2}-(n-k)}q^{k\choose 2}
S_{1,q/p}(n,k)\\
&=p^{{k\choose 2}-{n-k+1\choose 2}-(n-k)}S_{q/p}(n,k),
\end{align*}
we derive \eqref{eq:pqfrobenius} from \eqref{eq:qfrobenius} by
substituting $q\to q/p$ and $x\to xp^n$.

\small

\end{document}